\documentclass[reqno]{amsart}
\usepackage{pictex}
\usepackage{mathrsfs}
\usepackage{hyperref}
\usepackage{color}
\begin{document}

%
%

\def\labelenumi{(\theenumi)}

\newtheorem{thm}{Theorem}[section]
\newtheorem{lem}[thm]{Lemma}
\newtheorem{conj}[thm]{Conjecture}
\newtheorem{cor}[thm]{Corollary}
\newtheorem{add}[thm]{Addendum}
\newtheorem{prop}[thm]{Proposition}
\theoremstyle{definition}
\newtheorem{defn}[thm]{Definition}
\theoremstyle{remark}
\newtheorem{rmk}[thm]{Remark}
\newtheorem{example}[thm]{{\bf Example}}

\newcommand{\OmegaH}{\Omega/\langle H \rangle}
\newcommand{\hatOmegaHstar}{\hat \Omega/\langle H_{\ast}\rangle}
\newcommand{\SurfG}{\Sigma_g}
\newcommand{\TriangG}{T_g}
\newcommand{\TriangGOne}{T_{g,1}}
\newcommand{\ProjG}{\mathcal{P}_g}
\newcommand{\TeichG}{\mathcal{T}_g}
\newcommand{\CirclePackGTau}{\mathsf{CPS}_{g,\tau}}
\newcommand{\CrossRatio}{{\bf c}}
\newcommand{\CrossRatioGTau}{\mathcal{C}_{g,\tau}}
\newcommand{\CrossRatioOneTau}{\mathcal{C}_{1,\tau}}
\newcommand{\DeformGTau}{\mathcal{C}_{g,\tau}}
\newcommand{\Forget}{\mathit{forg}}
\newcommand{\Uniform}{\mathit{u}}
\newcommand{\Section}{\mathit{sect}}
\newcommand{\SLTwoC}{\mathrm{SL}(2,{\mathbb C})}
\newcommand{\SLTwoR}{\mathrm{SL}(2,{\mathbb R})}
\newcommand{\SUTwo}{\mathrm{SU}(2)}
\newcommand{\PSLTwoC}{\mathrm{PSL}(2,{\mathbb C})}
\newcommand{\GLTwoZ}{\mathrm{GL}(2,{\mathbb Z})}
\newcommand{\GLTwoC}{\mathrm{GL}(2,{\mathbb C})}
\newcommand{\PSLTwoR}{\mathrm{PSL}(2,{\mathbb R})}
\newcommand{\PGLTwoR}{\mathrm{PGL}(2,{\mathbb R})}
\newcommand{\GLTwoR}{\mathrm{GL}(2,{\mathbb R})}
\newcommand{\PSLTwoZ}{\mathrm{PSL}(2,{\mathbb Z})}
\newcommand{\SLTwoZ}{\mathrm{SL}(2,{\mathbb Z})}
\newcommand{\nnn}{\noindent}
\newcommand{\MCG}{{\mathcal {MCG}}}
\newcommand{\MMap}{{\bf \Phi}_{\mu}}
\newcommand{\HH}{{\mathbb H}^2}
\newcommand{\TT}{{\mathbb T}}
\newcommand{\X}{{\mathcal  X}}
\newcommand{\C}{{\mathscr C}}
\newcommand{\CC}{{\mathbb C}}
\newcommand{\RR}{{\mathbb R}}
\newcommand{\B}{{\mathcal  B}}
\newcommand{\G}{{\mathcal  G}}
\newcommand{\R}{{\mathcal  R}}
\newcommand{\Q}{{\mathbb Q}}
\newcommand{\ZZ}{{\mathbb Z}}
\newcommand{\PL}{{\mathscr {PL}}}
\newcommand{\GP}{{\mathcal {GP}}}
\newcommand{\GT}{{\mathcal {GT}}}
\newcommand{\GQ}{{\mathcal {GQ}}}
\newcommand{\EE}{{{\mathcal E}(\rho)}}
\newcommand{\HHH}{{\mathbb H}^3}

\def\square{\hfill${\vcenter{\vbox{\hrule height.4pt \hbox{\vrule width.4pt
height7pt \kern7pt \vrule width.4pt} \hrule height.4pt}}}$}

\newenvironment{pf}{\noindent {\sl Proof.}\quad}{\square \vskip 12pt}

\title[End Invariants for  $\SLTwoC$ characters]{ End Invariants for  $\SLTwoC$ characters of the one-holed torus}
\author{Ser Peow Tan, Yan Loi Wong, and Ying Zhang}
\address{Department of Mathematics \\ National University of Singapore \\
2 Science Drive 2 \\ Singapore 117543}
\email{mattansp@nus.edu.sg}%
\address{Department of Mathematics \\ National University of Singapore \\
2 Science Drive 2 \\ Singapore 117543}
\email{matwyl@nus.edu.sg}%
\address{Department of Mathematics \\Yangzhou University \\
Yangzhou, Jiangsu 225002 \\ China} \email{yingzhang@yzu.edu.cn}



\thanks{The authors are partially supported by the National University
of Singapore academic research grant R-146-000-056-112. The third
author is also partially supported by the National Key Basic
Research Fund (China) G1999075104.}

\vskip 20pt

%
%

\begin{abstract}
We define and study the set ${\mathcal E}(\rho)$ of end invariants
of a $\SLTwoC$ character $\rho$ of the one-holed torus $T$. We
show that the set $\EE$ is the entire projective lamination space
$\PL$ of $T$ if and only if (i) $\rho$ corresponds to the dihedral
representation, or (ii) $\rho$ is real and corresponds to a
$\SUTwo$ representation; and that otherwise, $\EE$ is closed and
has empty interior in $\PL$. For real characters $\rho$, we give a
complete classification of $\EE$, and show that $\EE$ has either
$0$, $1$ or infinitely many elements, and in the last case, $\EE$
is either a Cantor subset of ${\mathscr {PL}}$ or is ${\mathscr
{PL}}$ itself. We also give a similar classification for
``imaginary'' characters where the trace of the commutator is less
than 2. Finally, we show that for discrete characters (not
corresponding to dihedral or $\SUTwo$ representations), $\EE$ is a
Cantor subset of $\PL$ if it contains at least three elements.

\end{abstract}

\maketitle

\section{Introduction and statement of results}\label{s:intro}

Let $T$ be the one-holed torus, and $\pi$ its fundamental group
which is free on two generators $X, Y$. The $\SLTwoC$ character
variety of $T$ is the set $\X$  of equivalence classes of
representations $\rho:\pi \mapsto \SLTwoC$, where the equivalence
classes are obtained by taking the closure of the orbits under
conjugation by $\SLTwoC$. In this paper we define and study the
set of end invariants associated to the $\SLTwoC$ characters of
the one-holed torus. To simplify the exposition, by abuse of
notation, we use $\rho$ instead of $[\rho]$ to denote the
characters in ${\mathcal X}$ in the rest of the paper, there
should be no confusion, as we will be mostly  interested in the
trace function which is invariant under conjugation.

Let $\PL$ be the projective lamination space of $T$ and $\C
\subset \PL$ the set of (free homotopy classes of) essential
simple closed curves on $T$.

\begin{defn}\label{def:endinvariant}(End invariants)
An element $X \in \PL$ is  an end invariant of the character $\rho$
if there exists $K>0$ and a sequence of distinct elements $X_n \in
{\mathscr C}$ such that $X_n \rightarrow X$ and $|{\rm
tr}\,\rho(X_n)|<K$ for all $n$. 
\end{defn}

Denote by $\EE$ the set of end invariants of $\rho$. We show that
this is a closed subset of the projective lamination space $\PL$
of $T$ (Proposition \ref{prop:closed}); classify the characters
$\rho$ for which $\EE=\PL$ and show that otherwise $\EE$ has empty
interior in $\PL$ (Theorem \ref{thm:E=PL}); classify characters
$\rho$ for which $\EE=\emptyset$ (Theorem \ref{thm:E=emptyset});
and classify the possible structure of $\EE$ for reducible, real,
imaginary and discrete characters (Theorems
\ref{thm:Eforreducible}, \ref{thm:EforReal},
\ref{thm:EforImaginary} and \ref{thm:Efordiscrete}).

The set $\EE$ gives information about the dynamics of the action of
the mapping class group $\Gamma$ of $T$ on the character $\rho$, and
is closely related to the study of Kleinian groups, dynamical
systems, see for example \cite{goldman-stantchev} or
\cite{roberts1996physa}, and also certain problems in mathematical
physics, see \cite{roberts1996physa}.

\vskip 5pt

 The character variety $\X$ stratifies into relative
character varieties: for $\kappa \in \CC$, the $\kappa$-relative
character variety is the set of equivalence classes $\rho $ such
that
$${\rm tr}\,\rho(XYX^{-1}Y^{-1})=\kappa$$ for one (and hence any)
pair of generators $X,Y$ of $\pi$. Denote by ${\mathcal X}_{\kappa}$
the $\kappa$-relative character variety. By classical results of
Fricke (see for example \cite{goldman2004exposition} or
\cite{tan-wong-zhang2005gmm}), we have the following
identifications:
$${\mathcal X} \cong \CC^3, \qquad {\mathcal X}_{\kappa} \cong \{(x,y,z)\in \CC^3 ~|~
x^2+y^2+z^2-xyz-2=\kappa\},$$ with the identification given by
$$\iota:\rho \mapsto (x,y,z):=({\rm
tr}\rho(X), {\rm tr}\rho(Y), {\rm tr}\rho(XY)),$$ where $X, Y$ is a
fixed pair of generators of $\pi$. The topology on ${\mathcal X}$
and ${\mathcal X}_{\kappa}$ will be that induced by the above
identifications.

A character $\rho \in {\mathcal X}_{\kappa}$ such that $\iota(\rho)$
is a permutation of the triple $(0,0, \pm \sqrt{\kappa +2})$ is
called a \emph{dihedral character}, the image is generated by two
elliptics of order 2 and contains a cyclic subgroup of index 2.  A
character is said to be \emph{real} if $\iota(\rho)=(x,y,z) \in
\RR^3$, and \emph{imaginary} if two of the entries are purely
imaginary and the third real. We adopt the convention that for
imaginary characters, at least two of the entries are non-zero, so
that the dihedral characters are not imaginary; this will simplify
statements of results later. Real characters correspond to $\SLTwoR$
or $\SUTwo$ representations while imaginary characters correspond to
$\GLTwoR$ representations (see \cite{goldman1997am},
\cite{goldman2003gt} and \cite{goldman-stantchev}).

\vskip 5pt

The case $\kappa=-2$ corresponds to the so called type-preserving
representations and has been extensively studied in the context of
Kleinian groups, see for example \cite{bowditch1998plms},
\cite{keen-series1993t}, \cite{minsky1999am}.

\vskip 5pt

 The case $\kappa=2$ corresponds to the reducible
representations and is also somewhat special, it is convenient to
treat this case separately, see Theorem \ref{thm:Eforreducible}
and \S \ref{s:reducible}. Indeed, apart from Theorems
\ref{thm:E=PL}, \ref{thm:Eforreducible} and    \S
\ref{s:reducible}, the reader should consider all other results
and sections to exclude the reducible case, to avoid unnecessary
complications with exceptions caused by this case.


\vskip 5pt

 The mapping class group $\Gamma:=\pi_0({\rm
Homeo}^+(T)) \cong \SLTwoZ$ acts on ${\mathcal X}$ and ${\mathcal
X}_{\kappa}$ respectively, where we consider homeomorphisms of $T$
which fix pointwise a neighborhood of the boundary. We are
interested in the large scale behavior of this action, particularly
in the case where the action is not proper; this is reflected in the
set of end invariants of $\rho$. The dynamics of the action was
classified for real ${\rm SL}(2)$ characters by Goldman in
\cite{goldman2003gt}, and a partial classification was given for
imaginary characters by Stantchev in his PhD thesis in Maryland, see
Goldman and Stantchev \cite{goldman-stantchev}.

\vskip 5pt

 The set ${\mathscr C}$ of free homotopy classes of essential simple
closed curves in $T$ can be thought of as the set of vertices of the
pants graph ${{\mathscr C}(T)}$, where two vertices $X,Y \in
{\mathscr C}$ are connected by an edge if and only if $X$, $Y$ have
geometric intersection number one in $T$. ${{\mathscr C}(T)}$ can be
realized as the completion of the Farey triangulation $\mathcal F$
of the hyperbolic plane $\HH$. In this way, ${\mathscr C}$ is
naturally identified with $\hat \Q:=\Q \cup \{\infty\}$, and the
projective lamination space $\PL$ of $T$ is identified with the
projective real line $\hat \RR:=\RR \cup \{\infty\}$, the boundary
of the hyperbolic plane $\HH$. The mapping class group $\Gamma$ acts
on these sets and ${{\mathscr C}(T)}$ in a natural way, this action
is realized  via the isomorphism of $\Gamma$ with ${\rm SL}(2,
\mathbb Z)$, which acts on the upper half-plane as a model of $\HH$.

%
%
%
%

We now give the exact statements of our results. The first result
describes all characters $\rho$ for which $\EE =\PL$, and shows that
otherwise, $\EE$ has empty interior.

\begin{thm}\label{thm:E=PL}
The set of end invariants $\EE$ is equal to $\PL$ if and only if
{\rm(i)} $\rho$ is dihedral; or {\rm(ii)} $\rho$ corresponds to a
${\rm SU}(2)$ representation. Furthermore, if $\EE \neq \PL$, then
$\EE$ has empty interior in $\PL$.
\end{thm}

The above can be thought of as the opposite extreme of the following
theorem, characterizing the characters for which $\EE $ is empty,
which is a consequence of results in \cite{bowditch1998plms}
(Theorem 2), \cite{tan-wong-zhang2005gmm} (Theorem 2.3, Proposition
2.4) and \cite{tan-wong-zhang2004nsc} (Theorem 1.6); we will give a
sketch of the proof in \S \ref{s:proofofthmE=PL}.

\begin{thm}\label{thm:E=emptyset}{\rm(Bowditch, Tan-Wong-Zhang)}
The set of end invariants $\EE$ is empty if and only if $\rho$
satisfies

\begin{enumerate}
\item[(i)]  ${\rm tr}\rho(X) \not\in (-2,2)$ for all $X \in
{\mathscr C}$;

\item[(ii)] $|{\rm tr}\rho(X)| \le 2$ for only finitely many
{\rm(}possibly no{\rm)} $X \in {\mathscr C}$.
\end{enumerate}
\end{thm}

We call conditions (i) and (ii) in  Theorem \ref{thm:E=emptyset}
\emph{the extended BQ-conditions}.

\vskip 5pt

The reducible characters ($\kappa=2$) are somewhat special; the
following result classifies $\EE$ for such characters.

\begin{thm}\label{thm:Eforreducible}{\rm(End invariants for reducible
characters)} 
For $ \rho \in {\mathcal X}_2$,   $\EE=\{X_0\}$ or $\PL$. %
Furthermore, in the first case, if $X_0 \in {\mathscr C}$, then
${\rm tr}\, \rho(X_0) \in [-2,2]$ and ${\rm tr}\, \rho(X) \not\in
[-2,2]$ for all $X \in {\mathscr C}\setminus \{X_0\}$, while if $X_0
\not\in {\mathscr C}$, then ${\rm tr}\, \rho(X) \not\in [-2,2]$ for
all $X \in {\mathscr C}$; and in the second case, ${\rm tr}\,
\rho(X) \in [-2,2]$ for all $X \in {\mathscr C}$. %
\end{thm}

Note that in particular, $\EE$ is never empty in this case, so that
a reducible character never satisfies the extended BQ-conditions.

Denote by ${\mathcal X}^{\RR}$ and ${\mathcal X}_{\kappa}^{\RR}$
the real ${\rm SL}(2)$ character variety and relative character
varieties respectively.
We have the following classification of $\EE$ for $\rho \in
{\mathcal X}_{\kappa}^{\RR}$, together with the description of the
corresponding $\rho$; we exclude the case $\kappa=2$ which was
covered in the preceding theorem.

\begin{thm}\label{thm:EforReal}{\rm(End invariants for real
characters)} %
Suppose $\rho \in {\mathcal X}_{\kappa}^{\RR}$, with $\kappa \neq
2$. Then exactly one of the following must hold:
\begin{enumerate}
\item [(a)] $\EE =\emptyset$, and  $\rho$ satisfies the extended
BQ-conditions.

\item[(b)] $\EE =\{\hat  X\}$ where $\hat X \in {\mathscr C}$,
$\rho$ is a $\SLTwoR$ representation, ${\rm tr}\,\rho(\hat X) \in
(-2,2)$, and ${\rm tr}\,\rho( X) \not\in (-2,2)$ for all $X \in
{\mathscr C} \setminus \{ \hat X \}$.

\item[(c)] $\EE$ is a Cantor subset of $\PL$,  $\rho$ is a
$\SLTwoR$ representation,   ${\rm tr}\,\rho(X) \in (-2,2)$ for at
least two distinct $X \in {\mathscr C}$, and ${\rm tr}\,\rho(Y)
\not \in (-2,2) \cup \{\pm \sqrt{\kappa+2}\}$ for some element $Y
\in {\mathscr C}$.

\item[(d)] $\EE=\PL$, and $\rho$ satisfies the conditions of
Theorem \ref{thm:E=PL}, that is, $\rho$ is the dihedral
representation or a $\SUTwo$ representation.
\end{enumerate}

\noindent Furthermore, case {\rm (a)} occurs only when $ \kappa \in
(-\infty, 2)\cup [\,18, \infty)$; case {\rm (b)} when $ \kappa \in
[\,6, \infty)$; case {\rm (c)} when $ \kappa \in (2,\infty) $; and
case {\rm (d)} when $ \kappa \in [-2,2)\cup (2,\infty]$.
\end{thm}

For  $\rho\in {\mathcal X}^{\RR}$, Theorem \ref{thm:EforReal}
implies that if $\EE$ has more than one element, then $\EE$ is
either a Cantor set or all of $\PL$. Furthermore, if $\EE$ has only
one element $X$ (and $\kappa \neq 2$), then $X$ is rational, that
is, corresponds to a simple closed curve. These results are not true
for general complex characters, for example, punctured torus groups
with two geometrically infinite ends have two end invariants, and
those with one geometrically infinite end have an end invariant
which is irrational. Note also that condition (ii) of the extended
BQ-conditions follows from condition (i) in the case of real
characters with $\kappa \neq 2$; this will follow from the proof of
the theorem.

There is also an intriguing connection between the end invariants of
real characters  and the dynamical spectrum of the almost periodic
Schr\"odinger operator; see \cite{bowditch1998plms} and
\cite{roberts1996physa} for details. Roughly speaking, consider a
one (real) parameter family of characters  $\rho_E \in {\mathcal
X}_{\kappa}^{\RR}$ (parametrized by the energy $E \in {\mathbb R}$)
such that $\iota(\rho_E)=(2, E-V_0,E-V_1)$, where $V_0\neq V_1$ are
fixed constants and $\kappa=2+(V_0-V_1)^2$; and an irrational
element $\lambda \in \PL \setminus {\mathscr C}$. Then the set of
values of $E$ for which $\lambda \in {\mathcal E}(\rho_E)$
corresponds to the dynamical spectrum of a discrete almost periodic
Schr\"odinger operator, and the conjecture is that this set is
always a Cantor set.


\vskip 5pt

Denote by ${\mathcal X}^I$ and ${\mathcal X}_{\kappa}^I$ the
imaginary character variety and relative character varieties
respectively. Recall that dihedral characters are not in ${\mathcal
X}^I$ by our convention. Note that $\kappa \in \RR$ in this case. We
classify  $\EE$ for $\kappa<2$:

\begin{thm}\label{thm:EforImaginary}{\rm(End invariants for imaginary
characters)}
\begin{enumerate}
\item [(i)] $\kappa =-2:$ For $ \rho \in {\mathcal X}_{-2}^I$,
$\EE$ is either a Cantor subset of $\PL$, or consists of a single
element $X$ in ${\mathscr C}$. In the latter case, ${\rm
tr}\,\rho(X)=0$ and $\rho$ is equivalent under the action of the
modular group $\Gamma$ to a character corresponding to the triple
$(0, x, ix)$ where $x \in \RR$ satisfies $|x| \ge 2$.

\item[(ii)] $-14 <\kappa < 2:$ For $ \rho \in {\mathcal
X}_{\kappa}^I$, $\EE$ is either a Cantor subset of $\PL$, or
consists of a single element $X$ in ${\mathscr C}$.

\item[(iii)] $\kappa \le -14:$ For $ \rho \in {\mathcal
X}_{\kappa}^I$, $\EE$ is  a Cantor subset of $\PL$; consists of a
single element $X$ in ${\mathscr C}$; or is empty. %
\end{enumerate}
\end{thm}

We chose in the statement of Theorem \ref{thm:EforImaginary} above
to emphasize the case $\kappa=-2$ since the results are somewhat
sharper than for general $-14 \le \kappa <2$, and the case is itself
of independent interest.

\vskip 5pt

Our final result is for discrete characters. We say a character
$\rho$ is \emph{discrete} if the set of values  $\{{\rm tr}\,\rho(X)
\mid X \in {\mathscr C}\}$ is a discrete subset of $\CC$. Denote the
set of discrete characters by ${\mathcal X}^{\rm disc}$. We have:

\begin{thm}\label{thm:Efordiscrete}
For $ \rho \in {\mathcal X}^{\rm disc}$, if  $\EE$ has at least
three elements and $\EE \neq \PL$, then $\EE$ is  a Cantor set.
\end{thm}

The above can be rephrased as follows: For discrete characters not
corresponding to dihedral or $\SUTwo$ representations, $\EE$ is a
Cantor subset of $\PL$ if it contains at least three elements. Note
that a discrete character may also have 0, 1 or 2 elements in $\EE$.

The definition of an end invariant given generalizes the definition
of a (geometrically infinite, or degenerate) end for the case where
the representation is type-preserving ($\kappa =-2$), and discrete
and faithful, which is the subject of intensive study in the last
decade, especially in relation to Thurston's Ending Lamination
Conjecture, proven by Minsky (for the punctured torus) in
\cite{minsky1999am}. Our definition is motivated by that given by
Bowditch  in \cite{bowditch1998plms}, where end invariants were
defined for type-preserving but not necessarily discrete or faithful
representations. In fact, this work was very much inspired by
\cite{bowditch1998plms} and grew out of our attempt to study and
develop the thread of ideas presented in \S 5 of
\cite{bowditch1998plms}.

Note however that our definition differs slightly from that used in
\cite{bowditch1998plms}; the difference is that accidental
parabolics were (isolated) end invariants there, whereas they are
not in ours. This slight variation simplifies the statements of the
results above. It also allows us to state the following conjecture,
of which the preceding results can be regarded as supporting
evidence.

\begin{conj}\label{conj:cantorset}{\rm(The dendrite conjecture)}
Suppose that $\EE$ has more than two elements. Then either
$\EE=\PL$ or $\EE$ is a Cantor subset of $\PL$.
\end{conj}

The above conjecture is a refinement and generalization of the
suggestion by Bowditch in \cite{bowditch1998plms} that for a generic
$\rho \in {\mathcal X}_{-2}$ not satisfying the BQ-conditions, $\EE$
should be a Cantor set. The ``convex hull'' of $\EE$ is a subtree of
the dual tree $\Sigma$ of ${{\mathscr C}(T)}$; the above conjecture
says that this tree should look like a dendrite, in the sense that
if it has more than two ends, than there should be infinite
branching at any end and all ends are not isolated. The statement
would have been somewhat more complicated, with several exceptions,
if accidental parabolics are considered to be end invariants.

\vskip 5pt

Note that in the cases considered by Bowditch, he did not really
have to worry about the case $\EE=\PL$ since for $\kappa=-2$ this
only occurs for the trivial dihedral character corresponding to
the quaternionic representation $\rho$ with
$\iota([\rho])=(0,0,0)$.

\vskip 5pt

Bowditch did not produce examples of characters for which $\EE$ was
a Cantor set. Theorem \ref{thm:EforReal} produces many such examples
for real characters with $\kappa >2$ and Theorems
\ref{thm:EforImaginary} and \ref{thm:Efordiscrete} produces many
non-real examples for general $\kappa$, for example the character
$\rho \in {\mathcal X}_{-2}$ with $\iota(\rho)=(0,1,i)$ has $\EE$ a
Cantor set by either theorem. To the best of our knowledge, these
are the first examples for which $\EE$ is known to be a Cantor set.

\vskip 5pt

There is also a generalization of the Ending Lamination Conjecture
for $\SLTwoC$ characters (as Bowditch conjectured for the
$\kappa=-2$ case) which can be stated as follows:

\begin{conj}\label{conj:endlamination}
Suppose that $\rho, \rho' \in {\mathcal X}_{\kappa}$ are such that
$\EE={\mathcal E}(\rho')$, $\EE$ has at least two elements, and $\EE
\neq\PL$. Then $\rho=\rho'$.
\end{conj}

We end this introduction with a few words about the generalizations
to arbitrary surfaces. The definition of $\EE$ can be extended
without much difficulty. The case of the four-holed sphere is
similar and the techniques given here should give similar results in
that case, although the analysis is generally more difficult. In
other cases, $\PL$ is homeomorphic to the sphere $S^n$ for some $n
\ge 2$ and a possible generalization of Theorem \ref{thm:E=PL} is
that $\EE$ has either full measure or measure zero. A possible
generalization of Conjecture \ref{conj:cantorset} would be that
$\EE$ is perfect, if it contains more than two elements. However,
these are just speculations and we do not have any insights into
these more general cases.

The rest of the paper is organized as follows. In \S
\ref{ss:notation} we give the notation and basic definitions to be
used in the rest of the paper. In \S \ref{ss:basicresults} we state
three key lemmas used for the proofs of the theorems. In \S
\ref{s:proofofthmE=PL}, we prove Theorems \ref{thm:E=PL} and
\ref{thm:E=emptyset}. In \S \ref{s:reducible} we consider reducible
characters and prove Theorem \ref{thm:Eforreducible}. In \S
\ref{s:proofsforrealanddiscrete} we prove Theorems
\ref{thm:EforReal} and \ref{thm:Efordiscrete} and in \S
\ref{s:proofforimaginary}, we  prove Theorem
\ref{thm:EforImaginary}. Finally, in the Appendix, we give a brief
description of the $\tau$-reduction algorithm of Goldman-Stantchev
in \cite{goldman-stantchev} for the imaginary characters, which
generalizes that used by Bowditch in \cite{bowditch1998plms} and
which is used in a crucial way in the proof of Theorem
\ref{thm:EforImaginary}.

\vskip 6pt

\noindent {\it Acknowledgements.} Part of this work was carried out
while the first named author was visiting the University of
Maryland, College Park, the University of Warwick, and the Tokyo
Institute of Technology, he would like to thank his hosts Bill
Goldman, Caroline Series and Sadayoshi Kojima and these institutions
for their hospitality. He would also like to thank Rich Schwartz,
Rich Brown, George Stantchev, Javier Aramayona, John Parker, Juan
Souto, Brian Bowditch and especially Bill Goldman, Caroline Series,
Makoto Sakuma and Greg McShane for many stimulating and useful
conversations. He would also like to thank Shigeru Mizushima for
help with a computer program to help visualize the characters and
their end invariants.




\vskip 10pt

\section{Notation and definitions}\label{ss:notation}

\vskip 5pt

In this section we introduce the notation and definitions to be used
in the rest of this paper. As in the introduction, let $T$ denote
the one-holed torus, that is, a torus with an open disk removed. Its
fundamental group $\pi$ is freely generated by two elements $X,Y$
corresponding to two simple closed curves on $T$ with intersection
number one.

\subsection{The (relative) character variety $\X$
(resp. $\X_{\kappa})$} \label{sss:character variety} The $\SLTwoC$
character variety is the set
$${\mathcal X}={\rm Hom}(\pi, \SLTwoC)/\!/\,\SLTwoC$$
where the quotient is the geometric invariant theory quotient by the
conjugation action. By abuse of notation, denote by $\rho$ (instead
of $[\rho]$) the elements of ${\mathcal X}$; we call them
\emph{characters} of $T$. For $\kappa \in \CC$, the
$\kappa$-relative character variety is the subset
$$\X_{\kappa}=\{\rho \in \X ~|~ {\rm
tr}\,\rho(XYX^{-1}Y^{-1})=\kappa\},$$ where $X,Y$ are generators
of $\pi$. By results of Fricke, it does not matter which pair of
generators $X,Y$ are used to define $\kappa$. We have
\begin{equation}\label{eqn:kappavariety}
\X \cong \CC^3, \qquad {\mathcal X}_{\kappa} \cong \{(x,y,z)\in
\CC^3 ~|~ x^2+y^2+z^2-xyz-2=\kappa\};
\end{equation}
the identification is given by
$$\iota:\rho \mapsto (x,y,z):=({\rm
tr}\rho(X), {\rm tr}\rho(Y), {\rm tr}\rho(XY)),$$ where $X, Y$ is a
fixed pair of free generators of $\pi$. Conversely,
$\iota^{-1}(x,y,z)$ can be realized by the following representation
(see \cite{goldman2003gt}):
$$\rho(X)=A:=\left(%
\begin{array}{cc}
  x & 1 \\
  -1 & 0 \\
\end{array}%
\right), \qquad \rho(Y)=B:=\left(%
\begin{array}{cc}
  0 & -\zeta\\
  \zeta^{-1} & y \\
\end{array}%
\right),$$ where $\zeta+\zeta^{-1}=z$.

\subsection{Topology of the (relative) character variety}\label{sss:topologyofvariety} %
The topology on $\X$ and $\X_{\kappa}$ will be that induced by the
identifications defined in (\ref{eqn:kappavariety}) respectively.

\subsection{Action of  the mapping class group $\Gamma$}\label{sss:modularaction} %
The mapping class group
$$\Gamma:=\pi_0({\rm Homeo}^+(T))\cong {\rm SL}(2, \mathbb Z) $$
acts on $\pi$ and hence on ${\mathcal X}$; the action is given by
$$\phi(\rho)=\rho \circ \phi^{-1},$$ where $\phi \in \Gamma$ and $\rho \in {\mathcal X}$.
The action is not effective, the kernel is generated by the
elliptic involution corresponding to $-I \in \SLTwoZ$, so the
effective action is by $\SLTwoZ/\pm I=\PSLTwoZ$.
The quantity $\kappa={\rm tr}\,\rho(XYX^{-1}Y^{-1})$ is preserved
under the action of $\Gamma$ by results of Nielsen (see
\cite{goldman2003gt}); hence $\Gamma$ also acts on the relative
varieties $\X_{\kappa}$. With the identification of ${\mathcal X}$
and ${\mathcal X}_{\kappa}$ with the complex varieties in
(\ref{eqn:kappavariety}), the action of $\Gamma$ is realized via
polynomial maps on these varieties; it is generated by the cyclic
permutation
\begin{equation}\label{eqn:permutations}
 c:(x,y,z) \mapsto (z,x,y)
\end{equation}
and the involution
\begin{equation}\label{eqn:xyztoxyz'}
    s:(x,y,z) \mapsto (y,x,xy-z)
\end{equation}
corresponding to ${\small \Big(%
\begin{array}{cc}
  1 & -1 \\
  1 & 0 \\
\end{array}%
\Big), \Big(%
\begin{array}{cc}
  0 & 1 \\
  -1 & 0 \\
\end{array}%
\Big) }\in \PSLTwoZ$ respectively.

\subsection{Sign change automorphisms}\label{sss:signchange}
There is a $\ZZ/2 \times \ZZ/2$ action on ${\mathcal X}$ (resp.
${\mathcal X}_{\kappa}$) generated by simultaneously changing the
signs of two of the entries of $\iota(\rho)=(x,y,z)$. Two characters
are equivalent under this action if and only if they correspond to
lifts of the same representation of $\pi$ into $\PSLTwoC$. The large
scale behavior of the action of $\Gamma$ on ${\mathcal X}$ and the
end invariants of $\rho$ are not affected by these sign change
automorphisms; nonetheless, it will be convenient to use them for
some local trace reduction arguments later.

\subsection{The  pants graph ${\mathscr C}(T)$ and the
projective lamination space $\PL$} \label{sss:pantsgraph}

Let ${\mathscr C}$ denote the set of free homotopy classes of
essential (nontrivial, non-boundary) simple closed curves on $T$,
which can be regarded as the set of vertices of the pants graph
${{\mathscr C}(T)}$ of $T$, where two vertices are connected by an
edge if and only if the corresponding curves have geometric
intersection number one. The pants graph ${{\mathscr C}(T)}$ can be
concretely realized as the completion $\bar {\mathcal F}$ of the
Farey tessellation ${\mathcal F}$ of the upper half-plane $\HH$ in
$\HH \cup \hat \RR$ (see Figure \ref{fig:Farey}); recall that
${\mathcal F}$ is the tessellation of $\HH$ by ideal triangles where
the edges are the translates of the infinite geodesic $(0,\infty)$
by $\PSLTwoZ$. Note that ${{\mathscr C}(T)}$ is not locally finite,
and every vertex has infinite degree. The space $\PL$ of projective
laminations on $T$ can be regarded as the completion of ${\mathscr
C}$;  in this way, $\PL$ can be identified with $\hat \RR$, where
the irrational points of $\RR$ correspond to projective laminations
which are not closed, and ${\mathscr C}$ is identified with $\hat
\Q$, once we fix an identification of $X,Y$ and $XY$ with $0,
\infty$ and $1$, where $X, Y$ are a fixed pair of generators of
$\pi$. There is a natural topology on $\PL$ induced from the
topology of $\hat \RR$ which agrees with the usual topology on $\PL$
regarded as the completion of ${\mathscr C}$. There is also a
natural orientation induced on $\PL$ from this identification, where
we use the usual orientation of $\hat \RR$. It is also convenient to
use the conformal unit disk model of the hyperbolic plane to
visualize all this; in this way, ${{\mathscr C}(T)}$ is a
triangulation of the unit disk $D$, and $\PL$ is identified with the
unit circle $S^1 \cong \hat \RR$, where the anti-clockwise direction
is positive. We use upper case letters 
to denote elements of ${\mathscr C}$, and more generally, elements
of $\PL$. Occasionally, we will also use $\lambda$ to denote an
element of $\PL \setminus {\mathscr C}$.

\subsection{The dual trivalent tree $\Sigma$}\label{sss:dualtree} %
\begin{figure}
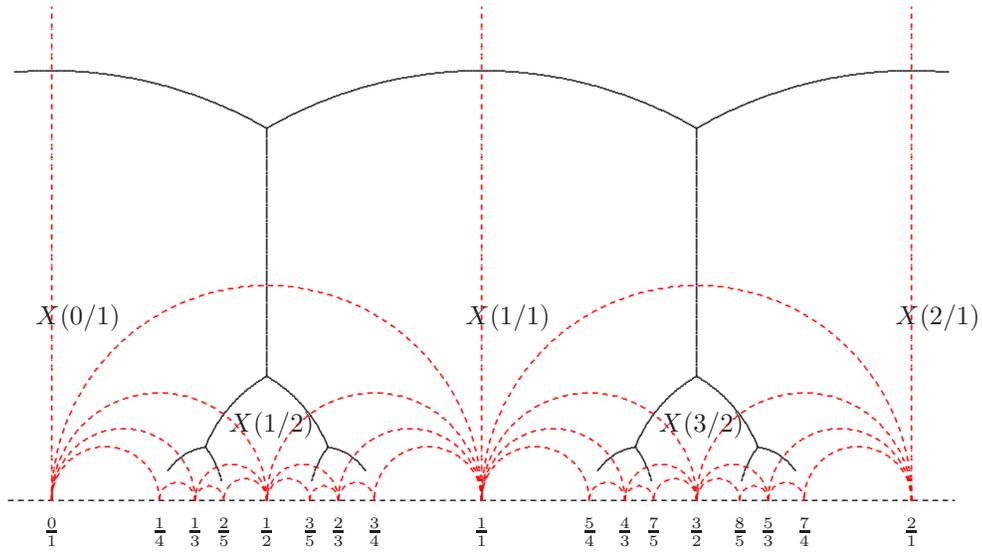


\begin{center}
\mbox{
\beginpicture
\setcoordinatesystem units <2.25in,2.25in>

\setplotarea x from -0.1 to 2.1, y from 0 to 1.2


\plot 0.5 0.86603 0.5 0.28868 /

\circulararc 38.21321 degrees from 0.5 0.28868 center at 0.66667 0

\circulararc -38.21321 degrees from 0.5 0.28868 center at 0.33333
0


 \circulararc 49.58256 degrees from 0.35714 0.12372  center at
0.375 0

 \circulararc -25.03966 degrees from 0.35714 0.12372 center
at 0.2 0

\circulararc 25.03966 degrees from 0.64286 0.12372 center at 0.8 0

\circulararc -49.58256 degrees from 0.64286 0.12372  center at
0.625 0


\plot 1.5 0.86603 1.5 0.28868 /

\circulararc 38.21321 degrees from 1.5 0.28868  center at 1.66667
0

\circulararc -38.21321 degrees from 1.5 0.28868 center at 1.33333
0

\circulararc 49.58256 degrees from 1.35714 0.12372  center at
1.375 0

\circulararc -25.03966 degrees from 1.35714 0.12372 center at 1.2
0

\circulararc 25.03966 degrees from 1.64286 0.12372 center at 1.8 0

\circulararc -49.58256 degrees from 1.64286 0.12372  center at
1.625 0


\circulararc 60 degrees from 1.5 0.86603  center at 1 0

\circulararc 35 degrees from 0.5 0.86603  center at 0 0

\circulararc -35 degrees from 1.5 0.86603  center at 2 0

\put {\mbox{\small $\frac{0}{1}$}} [cb] <0mm,-6mm> at 0 0

\put {\mbox{\small $\frac{1}{1}$}} [cb] <0mm,-6mm> at 1 0

\put {\mbox{\small $\frac{1}{2}$}} [cb] <0mm,-6mm> at 0.5 0

\put {\mbox{\small $\frac{1}{3}$}} [cb] <0mm,-6mm> at 0.33333 0

\put {\mbox{\small $\frac{2}{3}$}} [cb] <0mm,-6mm> at 0.66667 0

\put {\mbox{\small $\frac{1}{4}$}} [cb] <0mm,-6mm> at 0.25 0

\put {\mbox{\small $\frac{2}{5}$}} [cb] <0mm,-6mm> at 0.4 0

\put {\mbox{\small $\frac{3}{5}$}} [cb] <0mm,-6mm> at 0.6 0

\put {\mbox{\small $\frac{3}{4}$}} [cb] <0mm,-6mm> at 0.75 0


\put {\mbox{\small $\frac{2}{1}$}} [cb] <0mm,-6mm> at 2 0

\put {\mbox{\small $\frac{3}{2}$}} [cb] <0mm,-6mm> at 1.5 0

\put {\mbox{\small $\frac{4}{3}$}} [cb] <0mm,-6mm> at 1.33333 0

\put {\mbox{\small $\frac{5}{3}$}} [cb] <0mm,-6mm> at 1.66667 0

\put {\mbox{\small $\frac{5}{4}$}} [cb] <0mm,-6mm> at 1.25 0

\put {\mbox{\small $\frac{7}{5}$}} [cb] <0mm,-6mm> at 1.4 0

\put {\mbox{\small $\frac{8}{5}$}} [cb] <0mm,-6mm> at 1.6 0

\put {\mbox{\small $\frac{7}{4}$}} [cb] <0mm,-6mm> at 1.75 0

\put {\mbox{ $X(0/1)$}} [cb] <0mm,-6mm> at 0.05 0.5

\put {\mbox{ $X(1/1)$}} [cb] <0mm,-6mm> at 1.05 0.5

\put {\mbox{ $X(2/1)$}} [cb] <0mm,-6mm> at 2.05 0.5

\put {\mbox{ $X(1/2)$}} [cb] <0mm,-6mm> at 0.5 0.25

\put {\mbox{ $X(3/2)$}} [cb] <0mm,-6mm> at 1.5 0.25


\setdashes<1.9pt>

\plot -0.1 0 2.1 0 /

{\color{red}

\plot 0 0 0 1.15 / \plot 1 0 1 1.15 /

\circulararc 180 degrees from 1 0  center at 0.5 0

\circulararc 180 degrees from 0.5 0  center at 0.25 0

\circulararc 180 degrees from 1 0  center at 0.75 0

\circulararc 180 degrees from 0.33333 0  center at 0.16667 0

\circulararc 180 degrees from 0.5 0  center at 0.41667 0

\circulararc 180 degrees from 0.66667 0  center at 0.58333 0

\circulararc 180 degrees from 1 0  center at 0.83333 0

\circulararc 180 degrees from 0.25 0  center at 0.125 0

\circulararc 180 degrees from 0.33333 0  center at 0.29167 0

\circulararc 180 degrees from 0.4 0  center at 0.36667 0

\circulararc 180 degrees from 0.5 0  center at 0.45 0

\circulararc 180 degrees from 0.6 0  center at 0.55 0

\circulararc 180 degrees from 0.66667 0  center at 0.63333 0

\circulararc 180 degrees from 0.75 0  center at 0.70833 0

\circulararc 180 degrees from 1 0  center at 0.875 0

\plot 2 0 2 1.15 /

\circulararc 180 degrees from 2 0  center at 1.5 0

\circulararc 180 degrees from 1.5 0  center at 1.25 0

\circulararc 180 degrees from 2 0  center at 1.75 0

\circulararc 180 degrees from 1.33333 0  center at 1.16667 0

\circulararc 180 degrees from 1.5 0  center at 1.41667 0

\circulararc 180 degrees from 1.66667 0  center at 1.58333 0

\circulararc 180 degrees from 2 0  center at 1.83333 0

\circulararc 180 degrees from 1.25 0  center at 1.125 0

\circulararc 180 degrees from 1.33333 0  center at 1.29167 0

\circulararc 180 degrees from 1.4 0  center at 1.36667 0

\circulararc 180 degrees from 1.5 0  center at 1.45 0

\circulararc 180 degrees from 1.6 0  center at 1.55 0

\circulararc 180 degrees from 1.66667 0  center at 1.63333 0

\circulararc 180 degrees from 1.75 0  center at 1.70833 0

\circulararc 180 degrees from 2 0  center at 1.875 0

}

\endpicture
}\end{center}

\caption{Farey tessellation $\mathcal F$ and  the dual trivalent
tree $\Sigma$. }\label{fig:Farey}
\end{figure}

The dual graph of ${\mathcal F}$ (or ${\mathscr C}(T)$) is an
infinite trivalent tree $\Sigma$. Geometrically, we may choose the
vertices of $\Sigma$ as the incenters of the ideal triangles in the
tessellation ${\mathcal F}$ and the edges as the geodesic arcs
connecting the incenters of pairs of adjacent ideal triangles in
${\mathcal F}$, where two ideal triangles are said to be adjacent if
they share a common side. $\Sigma$ is trivalent since for each ideal
triangle in ${\mathcal F}$ there are exactly three others adjacent
to it; see Figure \ref{fig:Farey}. Note that $\Sigma$ is properly
imbedded in the hyperbolic plane and all its ends form the whole
ideal boundary of the hyperbolic plane. Denote by $V(\Sigma)$,
$E(\Sigma)$ the sets of vertices and edges of $\Sigma$ respectively,
and we use the notation $v$, $e$ to represent elements of
$V(\Sigma)$ and $E(\Sigma)$ respectively. A complementary region of
$\Sigma$ is the closure of a connected component of the complement
of $\Sigma$ in $\HH$; we denote the set of complementary regions of
$\Sigma$ by $\Omega(\Sigma)$.

\subsection{Generating pairs, triples and quadruples.}\label{sss:generatingpair} %
For $X,Y,Z,Z' \in {\mathscr C}$:
\begin{itemize}
\item The (unordered) pair $(X,Y)$, is a \emph{generating pair}
if $X$ and $Y$ are connected by an edge of ${{\mathscr C}(T)}$. We
say that $X$ and $Y$ are \emph{neighbors} (in ${{\mathscr C}(T)}$);

\item The (unordered) triple $(X,Y,Z)$ is a \emph{generating
triple} if $X,Y$ and $Z$ are the vertices of a triangle in
${{\mathscr C}(T)}$; and

\item $(X,Y;Z,Z')$ (where each of the first and second pair in the
quadruple is unordered) is a \emph{generating quadruple} if
$(X,Y,Z)$ and $(X,Y,Z')$ are generating triples (note that the pair
$Z,Z'$
is determined uniquely by $X,Y$). %
\end{itemize}

We shall denote the sets of generating pairs, triples and quadruples
by $\GP$, $\GT$ and $\GQ$ respectively.

\subsection{Correspondences between various sets of objects}
\label{sss:correspondence} There are natural correspondences between
the following sets, which are self-evident (see Figure \ref{fig:edge
oriented}):
\begin{eqnarray}\label{correspondences}
    V(\Sigma) &\longleftrightarrow &\GT~~;~~ v \mapsto (X,Y,Z) \\
    E(\Sigma) &\longleftrightarrow &\GP \longleftrightarrow \GQ~~;~~
    e\leftrightarrow (X,Y)\leftrightarrow (X,Y;Z,Z')\\
    \Omega(\Sigma) &\longleftrightarrow & {\mathscr
    C}\longleftrightarrow\hat \Q.
\end{eqnarray}

In the second correspondence, $(X,Y)$ defines an edge in $\mathcal
F$ which is dual to $e$. In the last correspondence, we will use the
same letters $X,Y,Z$ to denote the elements of all three sets,
namely, $\Omega(\Sigma)$, $\mathscr C$ and $\hat \Q$. Indeed, we
shall use the correspondence freely, so that the same symbol $X$ may
denote an element of $\Omega(\Sigma)$, $\mathscr C$  or $\hat \Q$;
it should be clear from the context which one we mean. We will also
use the symbols $X(p/q)$ to indicate that $X \in {\mathscr C}$ (or
$\Omega(\Sigma)$) corresponds to $p/q \in \hat \Q$.

\subsection{Directed edges of $\Sigma$} \label{sss:directededges}
\begin{figure}
\setlength{\unitlength}{1mm} 
\begin{picture}(60,40)
\thicklines \put(20,20){\vector(1,0){20}}
\put(20,20){\line(-1,1){10}} \put(20,20){\line(-1,-1){10}}
\put(40,20){\line(1,1){10}} \put(40,20){\line(1,-1){10}}
\put(30,30){$X$} \put(30,10){$Y$} \put(50,20){$Z'$}
\put(6,20){$Z$} \put(30,21){$\vec e$}
\end{picture}
\caption{The directed edge $\vec e \leftrightarrow (X,Y;Z
\rightarrow Z')$} \label{fig:edge oriented}
\end{figure}
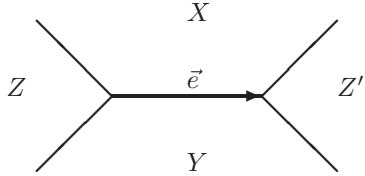

Let $\vec E(\Sigma)$ denote the set of directed edges of $\Sigma$.
Denote by $\vec e$ the elements of $\vec E(\Sigma)$, where the
direction of the arrow goes from the tail to the head. We use the
notation $\vec e \leftrightarrow (X,Y; Z \rightarrow Z')$ to
indicate that the directed edge $\vec e$ corresponds to the
generating quadruple $(X,Y;Z,Z')$ with the direction of the arrow
pointing from $Z$ towards $Z'$. We shall also use the orientation
convention that $X,Y,Z$ are in clockwise order as points on $\hat
\RR$ (as the boundary of $\HH$), so that if $\vec e \leftrightarrow
(X,Y; Z \rightarrow Z')$, then $-\vec e \leftrightarrow (Y,X; Z'
\rightarrow Z)$ where $-\vec e$ is the directed edge which is
directed in the opposite direction of $\vec e$ and has the same
underlying undirected edge $e$. In other words, directed edges
correspond to \emph{ordered} generating pairs, and \emph{ordered}
generating quadruples; see Figure \ref{fig:edge oriented}.


\subsection{Tri-coloring of ${\mathscr C}$, $\Omega(\Sigma)$ and $E(\Sigma)$}
\label{ss:tricolor}

The set ${\mathscr C}$(or $\Omega(\Sigma)$) can be naturally
partitioned into three equivalence classes, ${\R}$, $\G$ and $\B$
(we use the colors Red, Green and Blue to denote each class) in
such a way that any generating triple $(X,Y,Z)$ contains exactly
one element in each class. Similarly, $E(\Sigma)$ admits a
partition into three equivalence classes $E_r$, $E_g$ and $E_b$,
where if $e \leftrightarrow (X,Y)$, then $e \in E_r$ if $X,Y
\not\in \R$ and so on.

\subsection{Subsets of ${\mathscr C}$ and $\PL$}
\label{sss:subsetsofC}

Denote by $[X,Y]$ the set of points in $\PL$ in the closed interval
from $X$ to $Y$, going in the anti-clockwise/positive direction, so
for example $[X(0),X(1)]$ consists of all points in $\PL$
corresponding to real numbers in the interval $[0,1]$ and
$[X(1),X(0)]$ correspond to points in $\hat \RR$ outside the open
interval $(0,1)$. Denote similarly the half open and open subsets of
$\PL$ by $[X,Y)$, $(X,Y]$ and $(X,Y)$ (note that we use the same
notation for a generating pair but it should be clear from the
context if we mean a generating pair, or the corresponding open
interval in $\PL$). Denote by ${\mathscr C}_{[X,Y]}$ (and similarly
for the other intervals) the intersection ${\mathscr C} \cap [X,Y]$.
To each directed edge $\vec e \leftrightarrow (X,Y; Z \rightarrow
Z')$, we associate the subset of ${\mathscr C}$, called the tail of
$\vec e$ defined by  ${\rm tail}(\vec e):={\mathscr C}_{[X,Y]}$.
Note that this is the intersection of ${\mathscr C}$ with the closed
interval in $\PL$ with end points $X$ and $Y$ which contain $Z$.


\subsection{The Fricke trace map induced by a character $\rho$}\label{sss:tracemaps} %
A character $\rho \in {\mathcal X}_{\kappa}$ induces a map
$$\phi:=\phi_{\rho}:{\mathscr C} \rightarrow \CC$$
given by $\phi(X)={\rm tr}\,\rho(X)$. Equivalently,  $\phi$ can be
thought of as a map from $\Omega(\Sigma)$ to ${\mathbb C}$, from the
correspondence between $\C$ and $\Omega(\Sigma)$. We adopt the
convention that the corresponding lower case letters denote the
values of $\phi$, that is $\phi(X)=x$, $\phi(Y)=y$, $\phi(Z)=z$ and
so on. This will simplify notation considerably. This map is called
a \emph{Markoff map} by Bowditch in \cite{bowditch1998plms} when
$\kappa=-2$, and a \emph{generalized Markoff map} by the authors in
\cite{tan-wong-zhang2005gmm} for general $\kappa$ (where the map was
defined from $\Omega(\Sigma)$ to $\mathbb C$); we call this the
\emph{Fricke map} or \emph{Fricke trace map} in this paper. At any
rate, $\phi$ satisfies the following vertex and edge relations:
\begin{equation}\label{eqn:vertexrelation}
x^2+y^2+z^2-xyz-2=\kappa,
\end{equation}
for any generating triple $(X,Y,Z) \in \GT $; and
\begin{equation}\label{eqn:edgerelation}
    z+z'=xy,
\end{equation}
for any generating quadruple  $(X,Y;Z,Z') \in \GQ$ (recall we are
using the corresponding lower case letters to denote the values of
$\phi$).

The vertex relation (\ref{eqn:vertexrelation}) holds for every
generating triple if it holds for one by the edge relation
(\ref{eqn:edgerelation}), that is, it propagates along the tree
$\Sigma$. These relations arise from the corresponding Fricke
trace identities:
\begin{eqnarray}\label{eqn:2+tr[A,B]=}
({\rm tr}A)^2+({\rm tr}B)^2+({\rm tr}AB)^2-{\rm tr}A \,{\rm tr}B
\,{\rm tr}AB -2={\rm tr}[A,B],
\end{eqnarray}
\begin{eqnarray}\label{eqn:trAtrB=}
{\rm tr}AB+{\rm tr}AB^{-1}={\rm tr}A \,{\rm tr}B.
\end{eqnarray}

Note that the sign-change automorphism described in \S
\ref{sss:signchange} corresponds to changing the signs of
$\phi(X)$ for all $X$ in two of the classes in the tri-coloring of
$\C$ described in \S \ref{ss:tricolor},  while keeping the entries
in the third class fixed.

\subsection{The extended BQ-conditions.}\label{sss:BQConditions}
Let  $\phi$ be a Fricke trace map and $S \subseteq {\mathscr C}$.
We say that $\phi$ satisfies the extended BQ-conditions on $S$ if
\begin{enumerate}
\item[(i)]  $\phi(X) \not\in (-2,2)$ for all $X \in S$;

\item[(ii)] $|\phi(X)| \le 2$ for only finitely many (possibly no)
$X \in S$.
\end{enumerate}
If $S={\mathscr C}$, we say that $\phi$ satisfies the extended
BQ-conditions. If ${\mathscr C}=S_1 \cup \cdots \cup S_n$, then
clearly $\phi$ satisfies the extended BQ-conditions if and only if
it does on all the $S_i$, $i=1, \cdots, n$.


\subsection{Dihedral characters.}\label{sss:dihedral} Note that
$\iota^{-1}(0,0,\pm \sqrt{\kappa+2})$ is the representation
$$\rho(X)=A:=\left(%
\begin{array}{cc}
  0& 1 \\
  -1 & 0 \\
\end{array}%
\right), \qquad \rho(Y)=B:=\left(%
\begin{array}{cc}
  0 & -\zeta\\
  \zeta^{-1} & 0 \\
\end{array}%
\right),$$ where $\zeta+\zeta^{-1}=\pm \sqrt{\kappa+2}$. We call any
character for which two of the entries of $\iota(\rho)$ are zero
\emph{dihedral} characters; they are generated by
two order two elliptics $A$ and $B$ as above, and contain the cyclic subgroup {\small $\Big\langle \Big(%
\begin{array}{cc}
  \zeta & 0 \\
  0 & \zeta^{-1}\\
\end{array}%
\Big)\Big\rangle$} as a normal subgroup of index two. By
(\ref{eqn:edgerelation}), if $\rho$ is a dihedral character,
$\phi(X)$ takes only the values $0$, $\pm \sqrt{\kappa+2}$ for all
$X \in {\mathscr C}$.

\subsection{Real and imaginary characters}
\label{sss:realimaginarycharacters} A character $\rho$ is real if
$\iota(\rho) =(x,y,z) \in \RR^3$. Real characters correspond to
$\SLTwoR$ or $\SUTwo$ representations, see \cite{goldman1997am} or
\cite{goldman2003gt}. The real (relative) $\kappa$-character variety
can be identified with the real algebraic variety
\begin{equation}\label{eqn:realcharacters}
    {\mathcal X}_{\kappa}^{\RR}\cong \{(x,y,z) \in \RR^3 ~|~
    x^2+y^2+z^2-xyz-2=\kappa\}.
\end{equation}
A real character corresponds to a $\SUTwo$ representation if and
only if for any generating triple $(X,Y,Z)$, $-2 \le x, y, z \le 2$,
and $\kappa \le 2$ (see \cite{goldman2003gt}). Note that $ \kappa
\ge -2$ in this case. For an $\SLTwoR$ character which is not
dihedral, it is easy to see that for any $K>0$, there exists $X \in
{\mathscr C}$ with $|x|>K$ (see \cite{goldman2003gt}).

A character $\rho$ is imaginary if $\iota(\rho)$ is a triple with
two purely imaginary entries and one real entry, and at least two
entries are non-zero. If we use $ix$ and $iy$ for the imaginary
entries and $z$ for the real entry, where $x,y,z \in \RR$, then
the elements of ${\mathcal X}^I_{\kappa}$ can be identified with
the points on the real algebraic variety (taking away the points
on the coordinate axes since we are assuming the character is not
dihedral)
\begin{equation}\label{eqn:imaginarycharactervariety}
    \{(x,y,z) \in \RR^3
    ~|~-x^2-y^2+z^2+xyz-2=\kappa\}.
\end{equation}

\subsection{Reducible characters}\label{sss:reducible} A
character $\rho$ is reducible if $\kappa=2$ (see
\cite{goldman2003gt}). In this case, $\rho(X)$ and $\rho(Y)$ have a
common fixed point, which we may assume to be $\infty$, so that
$\rho$ is conjugate to an upper triangular representation. We can
replace $\rho$ with its semi-simplification which is a
representation by diagonal matrices with the same character, and
henceforth we shall do so. Hence the representations we shall
consider for reducible characters $\rho$ with $\iota(\rho)=(x,y,z)$
will be of the form
$$\rho(X)=\left(%
\begin{array}{cc}
  \xi & 0 \\
  0 & \xi^{-1} \\
\end{array}%
\right), \quad \rho(Y)=\left(%
\begin{array}{cc}
  \eta & 0 \\
  0 & \eta^{-1} \\
\end{array}%
\right), \quad
\rho(XY)=\left(%
\begin{array}{cc}
  \zeta & 0 \\
  0 & \zeta^{-1} \\
\end{array}%
\right),$$ where $\xi+\xi^{-1}=x,~~ \eta+\eta^{-1}=y,~~$ and $
\zeta+\zeta^{-1}=z$.

\subsection{Discrete characters}\label{sss:discretecharacters} A
character $\rho \in {\mathcal X}_{\kappa}$ is said to be discrete if
the set $\{\phi(X) ~|~ X \in {\mathscr C}\}$ is a discrete subset of
$\CC$. Note that this occurs for example if $\iota(\rho)=(x,y,z)$
and $x,y$ and $z$ all lie in a discrete ring in $\CC$, for example
$\ZZ$ or $\ZZ[i]$, or more generally the ring of integers $\mathcal
O_d$ of $\Q(\sqrt d)$, where $d<0$ is a square-free integer.

\subsection{The flow on $\Sigma$ associated to a character
$\rho$}\label{sss:flowofrho} Associated to a character $\rho$ is a
\emph{flow} on the tree $\Sigma$, that is, a map $f:=f_{\rho}:
E(\Sigma) \to \vec E(\Sigma)$ defined as follows: For
$e\leftrightarrow (X,Y;Z,Z')$, $f(e)=\vec e \leftrightarrow (X,Y; Z
\rightarrow Z')$ if $|z| \ge |z'|$, that is, the flow goes from the
larger absolute value of $\phi$ to the smaller one along $\Sigma$.
If $|z|=|z'|$, the flow is defined arbitrarily for that edge; this
ambiguity does not have any serious consequences on the subsequent
developments. A vertex $v \in V(\Sigma)$ is called a \emph{sink } if
$f(e)$ is directed towards $v$ for all the three edges $e \in
E(\Sigma)$ meeting at $v$. A finite subtree $\Sigma_0$ of $\Sigma$
is called an\emph{ attractor} for the flow if $f(e)$ is directed
towards $\Sigma_0$ for all $e \in \Sigma \setminus \Sigma_0$. Note
that a sink $v \leftrightarrow (X,Y,Z)$ may or may not be an
attractor for the flow; however, it is if in addition
$|x|,|y|,|z|>2$.

\subsection{End invariants}\label{sss:endinvariants} End
invariants for a character $\rho$ are defined as in Definition
\ref{def:endinvariant} in the introduction; the set of end
invariants is denoted by $\EE$. We shall see in \S
\ref{s:proofsforrealanddiscrete} that there are other equivalent
definitions which we will use. The relation between the end
invariants and (non-)properness of the action of $\Gamma$ on $\rho$
can be described as follows. $X \in \PL$ is an end invariant of
$\rho \in {\mathcal X}_{\kappa}$ if and only if there exists a
sequence $\{\theta_n\} \subset \Gamma \cong \SLTwoZ$ with $|{\rm
tr}~\theta_n|\rightarrow \infty$ such that $\theta_n(\rho)$ stays in
a fixed compact subset $S$ of ${\mathcal X}_{\kappa}$ (independent
of $\rho$), and the repelling fixed points $\mu_n^- \in \PL$ of
$\theta_n$ approach $X$ as $n \rightarrow \infty$.


\vskip 12pt

\section{Basic results}\label{ss:basicresults}

\vskip 3pt

Much of the ensuing discussion in this paper hinges on  three fairly
elementary but fundamental results on quasi-convexity (Lemma
\ref{lem:quasiconvexity}), escaping orbits (Lemma \ref{lem:escape})
and behavior of neighbors around $X$ (Lemma \ref{lem:neighbors})
which are easy to prove but play key roles in controlling the large
scale behavior of the action of  $\Gamma$ on the character $\rho$,
and in the proof of the theorems. They were first proved by Bowditch
in \cite{bowditch1998plms} for $\kappa=-2$; the generalization to
arbitrary $\kappa$ can be found in \cite{tan-wong-zhang2005gmm}.

For the rest of this section, we fix a $\rho \in \X_{\kappa}$ and
let the corresponding Fricke trace map be $\phi:=\phi_{\rho}$, which
we take to be a map from  $\Omega(\Sigma)$ to ${\mathbb C}$. We also
assume $\kappa \neq 2$. Recall that we adopt the convention
$\phi(X)=x$, $\phi(Y)=y$, and so on.

\subsection{Quasi-convexity: Connectedness of ${\mathscr
C}_{\phi}(K)$ for $K \ge 2$}\label{sss:quasiconvexity} %

We say that a subset $S \subset {\mathscr C}$ is connected if the
subgraph spanned by $S$ in the pants graph ${{\mathscr C}(T)}$ is
connected. For $K>0$, let ${\mathscr C}(K):={\mathscr
C}_{\phi}(K)=\{X \in {\mathscr C} ~|~ |\phi(X)| \le K\}$, and we
define $\Omega(K)$ similarly. We then have

\begin{lem}\label{lem:quasiconvexity}{\rm (Quasi-connectivity)}
For all $K \ge 2$, ${\mathscr C}(K)$ {\rm (}equivalently,
$\Omega(K)${\rm )} is connected.
\end{lem}

Lemma \ref{lem:quasiconvexity} can be deduced easily from the
following results.

\begin{prop}\label{prop:connectednessacrossedge}
Suppose that $K \ge 2$ and $(X,Y;Z,Z')$ is a generating quadruple
such that $Z,Z' \in {\mathscr C}(K)$. Then either $X$ or $Y$ {\rm
(}or both{\rm )} is in ${\mathscr C}(K)$.
\end{prop}

\begin{pf}
This follows directly from the edge relation
(\ref{eqn:edgerelation}).
\end{pf}

\begin{prop}\label{prop:twoarrowspointingout}
Suppose $X,Y,Z \in \Omega(\Sigma)$ meet at a vertex $v \in
V(\Sigma)$, and that the arrows on the edges $X \cap Y$ and $X
\cap Z$ arising from the flow $f_\rho$ both point away from $v$.
Then either $|x| \le 2$, or $y=z=0$.
\end{prop}

\begin{pf}
Let $Z'$ and $Y'$ be the regions opposite to $Z$ and $Y$
respectively, from the vertex $v$. By the assumption of the
direction of the arrow on the edge $X \cap Y$, we have $2|z| \ge
|z|+|z'| \ge |z+z'|=|xy|$. Similarly, $2|y| \ge |xz|$. Adding, we
get $2(|z|+|y|) \ge |x|(|y|+|z|)$ from which the conclusion follows.
\end{pf}

\noindent {\sl Sketch of proof of Lemma \ref{lem:quasiconvexity}}. \quad %
We prove connectedness of $\Omega(K)$. Suppose that this is not
connected. Take a minimal path in $\Sigma$ connecting two of the
components. If this path consists of only one edge in $\Sigma$, we
get a contradiction by Proposition
\ref{prop:connectednessacrossedge}. If it consists of more than one
edge, then by the construction, the flow at edges on the two ends of
the path point outwards, so that there is a vertex inside this path
where two of the arrows are pointing outwards, and we get a
contradiction by Proposition \ref{prop:twoarrowspointingout}. \qed

\subsection{Escaping orbits}

\begin{lem}\label{lem:escape}{\rm(Escaping orbits)}
Suppose that $\{\vec e_n\}$, $n \in \mathbb N$ is an infinite
directed path  in $\Sigma$ with the head of $\vec e_n$ equal to the
tail of $\vec e_{n+1}$, and such that $f_\rho(e_n)=\vec e_n$.
Furthermore, suppose that $\{\vec e_n\}$ does not limit to a
rational point of $\PL$. Then there exists infinitely many $X \in
{\mathscr C}(2)$ such that the path $\{\vec e_n\}$ intersects the
boundary of the corresponding complementary regions $X \in
\Omega(\Sigma)$.
\end{lem}

\begin{pf}
See \cite{bowditch1998plms} for the case $\kappa=-2$ and
\cite{tan-wong-zhang2005gmm} for the extension to the general case
where $\kappa \neq 2$.
\end{pf}

\vskip 6pt %
\subsection{Neighbors around $X$}~~ For each $X \in {\mathscr
C}$, let $Y_n$, $n \in \ZZ$ be the consecutive neighbors of $X$,
so that $(X,Y_n,Y_{n+1})\in \GT$ is a generating triple for all
$n$. For example, if $X$ corresponds to $\infty \in \hat \Q$, then
we can take $Y_n$ to correspond to $n \in {\mathbb Z}$. Let
$x=\lambda + \lambda^{-1}$ where $|\lambda| \ge 1$. Note that
$|\lambda|=1$ if and only if $x \in [-2,2] \subset {\mathbb R}$.
If $x=2$, then from the vertex relation (\ref{eqn:vertexrelation})
and edge relation (\ref{eqn:edgerelation}), $y_{n+1}=y_n \pm
\sqrt{\kappa-2}$, and $y_{n+1}-y_n=y_n-y_{n-1}$, hence the $\pm$
sign is constant in $n$. Similarly, if $x=-2$, then $y_{n+1}=-y_n
\pm \sqrt{\kappa-2}$, but this time, $y_{n+1}+y_n=-(y_n+y_{n-1})$,
hence the $\pm$ sign alternates in $n$. If $x=\pm
\sqrt{\kappa+2}$, then $y_{n+1}= \lambda^{\pm 1}y_n$ where the
$\pm$ sign is constant in $n$. If $x \notin \{\pm 2, \pm
\sqrt{\kappa+2}\}$ then there are (non-zero) constants $A,B \in
\mathbb C \backslash \{0\}$ with $AB=(x^2-\kappa-2)/(x^2-4)$ such
that $y_n=A {\lambda}^n + B {\lambda}^{-n}$. Hence we deduce that
the following holds. (This is Corollary 3.3 in
\cite{bowditch1998plms} in the case $\kappa=-2$.)

\begin{lem}\label{lem:neighbors}{\rm (Neighbors of $X$)} %
Suppose that $X\in {\mathscr C}$ has consecutive neighbors $Y_n$, $n
\in {\mathbb Z}$. Let $\rho \in {\mathcal X}_{\kappa}$, $\kappa \neq
2$, with corresponding Fricke trace map $\phi$.
\begin{enumerate}
\item [(a)] If $x \notin [-2,2] \cup \{\pm \sqrt{\kappa +2}\}$,
then $|y_n|$ grows exponentially as $n \rightarrow \infty$ and as
$n \rightarrow -\infty$.

\item [(b)] If $x \in (-2,2)$, then $|y_n|$ remains bounded. The
values of the neighbors $y_n$ are periodic if $x=2\cos \frac{p
}{q}\pi $ for some $p/q \in \mathbb Q$, and quasi-periodic
otherwise, that is, for any $n \in \mathbb Z$ and $\epsilon >0$,
there exists infinitely many indices $n_k$ with
$|y_n-y_{n_k}|<\epsilon$.

\item [(c)] If $x=2$, then either $y_n=y_0+n\sqrt{\kappa-2}$ for
all $n$, or $y_n=y_0-n\sqrt{\kappa-2}$ for all $n$. In particular,
since we assume that  $\kappa \neq 2$, $|y_n|$ grows linearly in
$|n|$.

\item [(d)] If $x=-2$, then either
$y_n=(-1)^{n}(y_0+n\sqrt{\kappa-2})$ for all $n$, or
$y_n=(-1)^{n}(y_0-n\sqrt{\kappa-2})$ for all $n$.

\item [(e)] If $x=\pm \sqrt{\kappa+2}$, then either
$y_n=\lambda^ny_0$ for all $n$, or $y_n=\lambda^{-n}y_0$ for all
$n$. In particular, if $x \not\in [-2,2]$ and $y_0 \neq 0$, then as
$n \rightarrow \infty$,  $|y_n| \rightarrow \infty $ {\rm (}or
$0${\rm )} and as $n \rightarrow -\infty$, $|y_n| \rightarrow 0$
{\rm (}resp. $\infty${\rm )}. %
\end{enumerate}
\end{lem}




\vskip 10pt

\section{Proofs of Theorems \ref{thm:E=PL} and \ref{thm:E=emptyset}}
\label{s:proofofthmE=PL}

\vskip 5pt

Fix $\rho \in {\mathcal X}_{\kappa}$, and again assume that $\kappa
\neq 2$. For a  pair of elements $X,Y \in {\mathscr C}$, recall that
$\rho$ satisfies the extended BQ-conditions on ${\mathscr
C}_{[X,Y]}$ if conditions (i) and (ii) of Theorem
\ref{thm:E=emptyset} are satisfied for all $Z \in {\mathscr
C}_{[X,Y]}$ (\S \ref{sss:BQConditions}). It was shown in
\cite{tan-wong-zhang2004nsc} that if this holds, then a version of
the McShane's identity holds, in particular, a certain series
converges which implies that  for all $K>0$, the set $\{Z \in
{\mathscr C}_{[X,Y]} ~|~ |z| \le K\}$ is finite. Hence, we have

\begin{prop}\label{prop:BQimpliesnotendinvariant}
Suppose that $\rho$ satisfies the extended BQ-conditions on
${\mathscr C}_{[X,Y]}$, then $\EE \cap {\mathscr
C}_{(X,Y)}=\emptyset$.
\end{prop}
We shall see later (Proposition \ref{prop:rationalendinvariants})
that if in addition, $x,y \neq \pm \sqrt{\kappa+2}$, then in fact
$\EE \cap {\mathscr C}_{[X,Y]}=\emptyset$.

Now suppose that $f_{\rho}(e)=\vec e$, with $\vec e \leftrightarrow
(X,Y ;Z\rightarrow Z')$ and $X,Y \not\in {\mathscr C}(2)$, that is,
$|x|,|y|>2$. Then it was shown in \cite{bowditch1998plms} (see also
\cite{tan-wong-zhang2005gmm}) that  for all edges $e'\in E(\Sigma)$
lying in the component of $\Sigma \setminus \{\vec e\}$ at the tail
of $\vec e$, $f_{\rho}(e')$ points towards $\vec e$. In particular,
for all $Z \in {\mathscr C}_{[X,Y]}$, $|z| \ge {\rm min}(|x|,|y|)$
and hence ${\rm tail}(\vec e)={\mathscr C}_{[X,Y]}$ satisfies the
extended BQ-conditions. Therefore we have

\begin{prop}\label{prop:tailofe}
Suppose that $f_{\rho}(e)=\vec e$, with $\vec e \leftrightarrow
(X,Y ;Z \rightarrow Z')$ and $X,Y \not\in {\mathscr C}(2)$. Then
$\EE \cap {\mathscr C}_{(X,Y)} =\emptyset$.
\end{prop}

Using the above result and Lemma \ref{lem:neighbors}, it is easy
to determine which elements of ${\mathscr C}$ are in $\EE$, and
also to show that $\EE$ is closed in $\PL$.

\begin{prop}\label{prop:rationalendinvariants}
Suppose $X \in {\mathscr C}$. Then $X \in \EE$ if and only if $x
\in (-2,2) \cup \{\pm \sqrt{\kappa+2}\}$.
\end{prop}

\begin{pf}
It is clear from parts (b) and (e) of Lemma \ref{lem:neighbors} that
$X \in \EE$ if  $x \in (-2,2) \cup \{\pm \sqrt{\kappa+2}\}$. Now
suppose that $x \not\in (-2,2) \cup \{\pm \sqrt{\kappa+2}\}$. Then
by parts (a), (c) and (d) of Lemma \ref{lem:neighbors}, for all
$K>0$, there exists $N \in {\mathbb N}$ such that $|y_n|>K$ for all
$n>N$ and $n<-N$. Hence, for all $K>0$, by the remark preceding
Proposition \ref{prop:tailofe}, there exists a neighborhood $N(X)$
of $X$ (depending on $K$) such that $|z|>K$ for all $Z \in
N(X)\setminus \{X\}$. We conclude that $X \not\in \EE$.
\end{pf}

\begin{prop}\label{prop:closed}
For any character $\rho \in {\mathcal X}$, the set of end invariants
$\EE$ is closed in $\PL$.
\end{prop}

\begin{pf}
We show that $\PL \setminus \EE$ is open. Suppose first that $X \in
{\mathscr C}$ is not an end invariant. Then by Proposition
\ref{prop:rationalendinvariants} and Lemma \ref{lem:neighbors}, if
$Y_n$ are the neighbors of $X$, $|y_n| \rightarrow \infty$ as $n
\rightarrow \pm \infty$. Hence, by Propositions \ref{prop:tailofe},
\ref{prop:rationalendinvariants}, there exists an open neighborhood
$U$ of $X$ such that $U \cap \EE =\emptyset$. Now suppose $\lambda
\in \PL \setminus {\mathscr C}$ is not an end invariant. Let
$\{e_n\}$ be a path in $\Sigma$ limiting to $\lambda$ and let
$(X_n,Y_n)\in \GP$ be the generating pair corresponding to the edge
$e_n$, oriented so that $\lambda \in [X_n,Y_n]$. Since $\lambda
\not\in \EE$, we have $|x_n|, |y_n| \rightarrow \infty$; so we may
assume say that $|x_1|=K>2$, and both $|x_N|$ and $|y_N|>K$ for some
$N \in \mathbb N$. By Lemma \ref{lem:quasiconvexity}, $|z|>K$ for
all $Z \in [X_N, Y_N]$, so that ${\mathscr C}_{[X_N, Y_N]}$
satisfies the extended BQ-conditions. Hence, $(X_N, Y_N) \cap \EE
=\emptyset $ by Proposition \ref{prop:BQimpliesnotendinvariant}. We
conclude that $\PL \setminus \EE$ is open.
\end{pf}

\noindent {\sl Proof of Theorem \ref{thm:E=PL}}. \quad If $\rho$ is
dihedral, then $x \in \{0, \pm \sqrt{\kappa +2}\}$ for all $X \in
{\mathscr C}$, and if $\rho$ is a ${\rm SU}(2)$ character, then $x
\in [-2,2]$ for all $X \in {\mathscr C}$. In either case, it is
clear that $\EE=\PL$ since $|x|$ is bounded for all $X \in {\mathscr
C}$. Now suppose that $\EE$ contains an open interval ${\mathcal I}
\subset \PL$. Then there exists a generating triple $(X,Y,Z) \in
\GT$ with $X,Y,Z \in {\mathcal I}$. By Proposition
\ref{prop:rationalendinvariants}, $x,y,z \in (-2,2) \cup \{\pm
\sqrt{\kappa+2}\}$. Suppose say that $x =\pm \sqrt{\kappa+2}$, where
$\sqrt{\kappa+2} \not\in[-2,2]$. Then either $y=z=0$ in which case
$\rho$ is a dihedral character, or $y,z \neq 0$, in which case by
part (e) of Lemma \ref{lem:neighbors} and Proposition
\ref{prop:rationalendinvariants}, there exists a neighbor $Y_n$ of
$X$ such that $Y_n \in {\mathcal I}$ but $Y_n \not\in \EE$, which
gives a contradiction. Hence, we may suppose that $x,y,z \neq \pm
\sqrt{\kappa+2}$, in which case $x,y,z \in(-2,2)$ (recall $\kappa
\neq 2$). If $\rho$ is a ${\rm SU}(2)$ character, $\EE=\PL$ and we
are done; otherwise, $\rho$ is a $\SLTwoR$ character and there
exists $W \in {\mathscr C}$ with $|w|>2$, $w \neq \pm \sqrt{\kappa
+2}$ (see for example \cite{goldman2003gt}) and so $W \not\in \EE$.
Since $x \in (-2,2)$, by Lemma \ref{lem:neighbors}(b), the
neighboring values around $X$ are either periodic or quasi-periodic;
it follows that there exists $W_n$ arbitrarily close to $X$ (hence
in ${\mathcal I}$) with $W_n \not\in \EE$. The contradiction
completes the proof. \qed

\vskip 10pt

\noindent {\sl Proof of Theorem \ref{thm:E=emptyset}}. \quad Suppose
that $\rho$ satisfies the extended BQ-conditions. It was shown in
\cite{tan-wong-zhang2004nsc} that in this case, there exists a
finite subtree of $\Sigma$ which is an  attractor for the flow
$f_{\rho}$ associated to $\rho$, and that the generalized McShane's
identity holds. This  implies that for any $K>0$, the set $\{X \in
{\mathscr C}~|~ |x| \le K\}$ is finite. It follows that
$\EE=\emptyset$. Note that Lemmas \ref{lem:quasiconvexity} and
\ref{lem:escape} played essential roles in the proof, in particular,
if the extended BQ-conditions are satisfied, there cannot be an
escaping orbit in the sense of Lemma \ref{lem:escape}, which is a
crucial step towards showing the existence of the attractor.

Conversely, if $\rho$ does not satisfy the extended BQ-conditions,
then either there exists some $X \in \C$ with $x \in (-2,2)$ or the
set $\C(2)=\{X \in \C ~|~ |x| \le 2\}$ is infinite. In the first
case, $X \in \EE$, and in the second case, $\C(2)$ has an
accumulation point in $\PL$ which lies in $\EE$. \qed


\vskip 15pt
\section{Reducible characters}\label{s:reducible}
\vskip 5pt

We consider the reducible characters in this section and prove
Theorem \ref{thm:Eforreducible}. As pointed out in \S
\ref{sss:reducible}, we may use representations into diagonal
matrices of $\SLTwoC$ to represent the reducible characters. Note
that $\rho(XYX^{-1}Y^{-1})=I$ for all generating pairs $X,Y$ of
$\pi$, where $I$ is the identity matrix. Hence we may think of
$\rho$ as corresponding to a  representation of the fundamental
group of the torus $\mathbb T$ (without boundary) and use the
homology classes $[X^mY^n]$ where $m,n$ are relatively prime and
$n \ge 0$ to represent the elements of ${\mathscr C}$. Let
$\pi({\mathbb T})=\langle X,Y ~|~ XYX^{-1}Y^{-1}=I \rangle $,
then $$\rho(X)=\left(%
\begin{array}{cc}
  \alpha^{1/2} & 0 \\
  0 & \alpha^{-1/2} \\
\end{array}%
\right), \qquad \rho(Y)=\left(%
\begin{array}{cc}
  \beta^{1/2} & 0 \\
  0 & \beta^{-1/2} \\
\end{array}%
\right),$$ where $\alpha, \beta \in \CC$, and $x={\rm tr}\,\rho(X)=
\alpha ^{1/2}+\alpha ^{-1/2}$, $y={\rm tr}\,\rho(Y)= \beta
^{1/2}+\beta ^{-1/2}$. There are two cases to consider, (i) $p\cdot
\log |\alpha |+q\cdot \log |\beta |=0$ for some integers $p,q$ (not
all zero); and (ii) otherwise. In the first case, we have
$|\alpha^p\cdot \beta^q|=1$. By using the action of $\Gamma$, we
might as well assume that in fact $p=1$ and $q=0$, so that
$|\alpha|=1$, and $x\in [-2,2]$. If $|\beta|=1$ as well, then $y\in
[-2,2]$ and similarly $z={\rm tr}\,\rho(XY) \in [-2,2]$. Since
$\kappa=2$, $\rho$ corresponds to a ${\rm SU}(2)$ representation and
$\EE=\PL$ in this case. If $|\beta| \neq 1$, then $|z| \not\in
[-2,2]$ for all other $Z \in \C \setminus \{X\}$, since the exponent
of $Y$ is non-zero for all the other elements. Furthermore, since
$|\alpha ^n \beta|=|\beta|$ is bounded, $X \in \EE$. For any other
element $Z \in \PL$, if $Z_n$ is any sequence of distinct elements
approaching $Z$, then the exponent of $Y$ in the homology class of
$Z_n$ approaches $\pm \infty$, so that $\log |Z_n| \to \pm \infty$.
Hence $Z \not\in \EE$.  It follows that in this case, $\EE=\{X\}$.
Now we consider case (ii), where $\frac{\log |\alpha |}{\log
|\beta|}=\mu$ is an irrational number. There exists a sequence of
rationals ${p_n}/{q_n} \rightarrow \mu$, namely, the convergents of
$\mu$ in the continued fraction expansion of $\mu$ such that
$p_n\cdot \log |\alpha |+q_n\cdot \log |\beta|$ is bounded. Hence
$X(\mu) \in \PL$ lies in $\EE$. For any other $Z \in \PL$, if
$X({p_n}/{q_n})$ is any sequence of elements approaching $Z$, then
$p_n\cdot \log |\alpha|+q_n\cdot \log |\beta|$ is unbounded and
hence $Z \not\in \EE$. It follows that $\EE=\{X(\mu)\}$ in this
case. This completes the proof of Theorem \ref{thm:Eforreducible}.



\vskip 15pt

\section{Proofs of Theorems \ref{thm:EforReal} and
\ref{thm:Efordiscrete}}\label{s:proofsforrealanddiscrete}

\vskip 5pt

We first give an equivalent definition for the set $\EE$ in terms of
the (projectivized) ends of a certain subtree ${\mathcal H}$ of
$\Sigma$; this is the variation of the original definition given by
Bowditch in \cite{bowditch1998plms}.

Suppose $\rho \in \X$  does not satisfy the extended BQ-conditions.
We define a subtree ${\mathcal H}:={\mathcal H}_{\rho}$ of $\Sigma$
as follows: Suppose that $e \in E(\Sigma)$ corresponds to the
generating pair $(X,Y)$. Then $e \subset {\mathcal H}$ if both
${\mathscr C}_{[X,Y]}$ and ${\mathscr C}_{[Y,X]}$ does \emph{not}
satisfy the extended BQ-conditions.

For edges $e \in E(\Sigma)\setminus {\mathcal H}$, exactly one of
the sets ${\mathscr C}_{[X,Y]}$, ${\mathscr C}_{[Y,X]}$ does not
satisfy the extended BQ-conditions (since $\rho$ doesn't); we define
a direction $d(e) \in \vec E(\Sigma)$ so that tail of $d(e)$
satisfies the extended BQ-conditions. Hence, $d(e)$ points towards
${\mathcal H}$ for all $e \not\in {\mathcal H}$. Note that $d(e)$ is
not to be confused with the flow $f(e)$ defined in \S
\ref{sss:flowofrho}. It is easy to see from the definition that any
vertex $v \in V(\Sigma)$ has $0$, $2$ or $3$ edges of ${\mathcal H}$
adjacent to it.  A vertex $v \in {\mathcal H}$ is called a
\emph{node} of ${\mathcal H}$ if there are three edges of ${\mathcal
H}$ adjacent to it. Furthermore, if $v$ is not a vertex of
${\mathcal H}$, then for the three edges incident at $v$, $d(e)$
points towards $v$ for two of the edges and away from $v$ for the
third. If $v$ has 2 edges of ${\mathcal H}$ adjacent to it, then
$d(e)$ points towards $v$ for the remaining edge.

There is a natural projection map from the ends of $\Sigma$ to $\PL$
which is one to one onto irrational points and two to one to
rational points. If ${\mathcal H}$ is empty then $d(e)$ all point
towards a unique $\lambda \in \PL$ and $\EE=\{\lambda\}$; otherwise,
$\EE$ is the image of the ends of ${\mathcal H}$ under the
projection map to $\PL$.

For an irrational $\lambda \in \EE$, the following result states
that we can take $K=2$ in the definition of $\EE$.

\begin{prop}\label{prop:irrationalend}
Suppose that $\lambda \in \EE$ and $\lambda \not\in {\mathscr C}$.
Then there exists a sequence of distinct elements $Z_n \in
{\mathscr C}$ such that $Z_n \rightarrow \lambda$ and $|z_n| \le
2$ for all $n$.
\end{prop}

\begin{pf}
If ${\mathcal H}$ is empty, then $d(e)$ points towards $\lambda$ for
all $e \in E(\Sigma)$. Choose any path $\{ e_n\}$ limiting at
$\lambda$ and let $(X_n,Y_n)$ be the generating pair corresponding
to $e_n$, ordered so that $\lambda \in [X_n,Y_n]$. Then $\rho$ does
not satisfy the extended BQ-conditions on ${\mathscr
C}_{[X_n,Y_n]}$, and hence, there exists $Z_n \in {\mathscr
C}_{[X_n,Y_n]}$ such that $|z_n|\le 2$. Since $X_n, Y_n \rightarrow
\lambda$, by passing to a subsequence if necessary, we obtain a
sequence of distinct $Z_n \rightarrow \lambda$ with $|z_n| \leq 2$.
The same argument works if $\lambda$ is the end of a path in
${\mathcal H}$: we just choose a path in $\mathcal H$ ending at
$\lambda$.
\end{pf}

\noindent {\sl Proof of Theorem \ref{thm:EforReal}}. \quad Let $\rho
\in {\mathcal X}_{\kappa}^{\mathbb R}$ with $\phi$ the corresponding
Fricke trace map,  and suppose that $\kappa \neq 2$. We first prove
the following claim:

\vskip 3pt

\noindent \emph{Claim}. \quad If $|\EE|>1$, then $|\EE \cap
\C|=\infty$, and if $\EE=\{X\}$, then $X \in \C$. Furthermore, if
$|\EE|>1$, then $\EE$ is perfect (that is, every element of $\EE$ is
an accumulation point of $\EE$).

\vskip 3pt

First note that $x=\phi(X) \in \mathbb R$ for all $X \in \C$ by
the edge relation (\ref{eqn:edgerelation}). Hence by Proposition
\ref{prop:rationalendinvariants}, $X \in \EE$ if $|x|<2$. If there
exists two distinct $X, Y \in \C$ with $|x|, |y| <2$, then  by
looking at the behavior of $\phi$ about $X$ and $Y$ respectively
and applying Lemma \ref{lem:neighbors} (b), we see that there are
infinitely many $Z \in \C \cap \EE$ and both $X$ and $Y$ are
accumulation points in $\EE$. Similarly, if there exists $X \in \C
\cap \EE$ with $x=\pm \sqrt{\kappa+2}$, and $\sqrt{\kappa+2}>2$,
then by Lemma \ref{lem:neighbors} (e), there are infinitely many
$Z \in \C \cap \EE$ and $X$ is an accumulation point of $\EE$ (in
fact, it is a one sided limit point of $\EE$ if $\rho$ is not
dihedral). Now suppose that $\lambda \in \EE$ with $\lambda
\not\in \C$. By Proposition \ref{prop:irrationalend},
$|\C(2)|=\infty$ and since $\C(2)$ is connected, we can find a
sequence $\{X_n\}\subset \C(2)$ such that $X_n \rightarrow
\lambda$ and $(X_n,X_{n+1})$ is a generating pair for all $n$.  If
there are infinitely many $X_k$ in the sequence such that
$|x_k|<2$, we are done. Otherwise, re-indexing if necessary, we
may suppose that $x_n=\pm 2$ for all $n \in \mathbb N$.  Note that
$X_{k-1}$ and $X_{k+1}$ are neighbors of $X_k$. If $(X_{k-1}, X_k,
X_{k+1})$ is a generating triple, then using the sign-change
automorphism (\S \ref{sss:signchange}), we may assume that
$(x_{k-1},x_{k},x_{k+1})=(2,2,2)$ or $(-2,-2,-2)$. In the first
case, $\kappa=2$ which contradicts our assumption. The second case
corresponds to the holonomy representation of the thrice punctured
sphere and $\phi$ satisfies the extended BQ-conditions, so
$\EE=\emptyset$, again a contradiction. Hence, $(X_{k-1}, X_k,
X_{k+1})$ is not a generating triple and there exists $Z_k$, a
neighbor of $X_k$ lying between $X_{k-1}$ and $X_{k+1}$. Again, by
the sign change automorphism, we may assume that $x_k=2$. Now by
Lemma \ref{lem:neighbors} (c), $x_{k-1}$ and $x_{k+1}$ have
opposite signs, say $x_{k-1}=-2$ and $x_{k+1}=2$, which forces
$|z_k|<2$. This produces a sequence of elements $Z_k \in \EE
\cap\, \C$ approaching $\lambda$ which completes the proof of the
claim.

The classification of the types for $\EE$ according to (a), (b), (c)
or (d) of the theorem now follows from the above claim and Theorems
\ref{thm:E=PL} and \ref{thm:E=emptyset}. It remains to show that
each of these cases occur for the values of $\kappa$ stated. We use
here Goldman's main Theorem in \cite{goldman2003gt} which classifies
the action of $\Gamma$ on ${\mathcal X}_{\kappa}^{\mathbb R}$. The
fact that case (a) occurs if and only if  $\kappa <2$ or $\ge 18$
follows from the fact that $\Gamma$ acts properly only on characters
in these ranges of $\kappa$. The ranges for case (c) and (d) also
follow from Goldman's result, Theorem \ref{thm:E=PL}, the above
argument, and the fact that the dihedral character is real for
$\kappa >-2$. It remains to show that case (b) occurs if and only if
$\kappa \ge 6$. It is easy to see that $\EE=\{X\}$ if
$\iota(\rho)=(0,a,a)$, where $a\ge 2$. Since in this case $\kappa=
2a^2-2 \ge 6$, we have case (b) occurs for all $\kappa \ge 6$.
Conversely, we show that if case (b) occurs, that is, if
$\EE=\{X\}$, then $\kappa \ge 6$. Let $Y_n$, where $n \in \mathbb
Z$, be the neighbors of $X$. Using the sign change automorphism, we
may assume that $x \ge 0$, so that $x \in [\,0,2)$,
$x=e^{i\theta}+e^{-i \theta}$ where $0<\theta\le \pi/2$, and there
exist $A,B \in \mathbb C$ with $AB={(x^2-\kappa-2)}/{(x^2-4)}$ such
that $y_n=Ae^{in\theta}+Be^{-in\theta}$. Since $y_n \in \mathbb R$
for all $n$, we have $B=\bar{A}$ so that $y_n=2\,{\rm Re}
(Ae^{in\theta})$. Now since $|y_n| \ge 2$ for all $n$, $\theta$ must
be a rational multiple of $\pi$, that is,  $\theta=p \pi/q$ for some
rational $p/q$, and the set $\{Ae^{i n \theta}\}$ is the set of
vertices of a regular $q$-gon in the complex plane centered at the
origin with $|2\,{\rm Re}(Ae^{i n \theta})|\ge 2$ for all $n$.
Re-indexing if necessary, we may assume that $y_0=2\,{\rm Re}(A)>0$;
hence $y_0 \ge 2$, and $y_1=2\,{\rm Re}(Ae^{i\theta})<0$, which
implies $y_1\le -2$. It follows that
$\kappa=x^2+y_0^2+y_1^2-xy_0y_1-2 \ge 8-2=6$ which completes the
proof of Theorem \ref{thm:EforReal}. \qed

\vskip 10pt


\noindent {\sl Proof of Theorem \ref{thm:Efordiscrete}}. \quad
Suppose that $\rho$ is discrete, and $\EE$ has at least three
elements but is not equal to $\PL$. By Theorem \ref{thm:E=PL}, it
suffices to prove that $\EE$ is perfect, that is every $X \in \EE$
is the limit of distinct $Z_n \in \EE$. Note that in this case
${\mathcal H}$ has at least three distinct ends and at least one
node.

{\it Case 1.  $X \in \C \cap \EE$}. \, By Proposition
\ref{prop:rationalendinvariants} and the discreteness of $\rho$, $x
=2\cos \frac{p}{q}\pi$ for some $\frac{p}{q}\in \mathbb Q$, and by
the periodicity of $\phi$ about $X$ ($\rho$ is stabilized by a
reducible element of the mapping class group $\Gamma$ fixing $X$)
and the fact that ${\mathcal H}$ has at least 3 ends, $X$ is a limit
of distinct $Z_n \in \EE$.

{\it Case 2:  $X \not\in \C$}. \, Choose a node $v$ of ${\mathcal
H}$ (which exists by the assumption) and the path $\{ e_n\}$ from
$v$ to $X$ and let $(X_n,Y_n)$ be the generating pair corresponding
to $e_n$, ordered so that $X \in [X_n,Y_n]$. By the connectedness of
$\C (2)$ (Lemma \ref{lem:quasiconvexity}) and the vertex relation
(\ref{eqn:vertexrelation}), it is easy to see that there exists a
universal constant $K>0$ depending only on $\kappa$ such that
$|x_n|, |y_n| <K$. Since the set of values are discrete, by passing
to a subsequence, we obtain a sequence of generating pairs
$(X_n,Y_n)$ approaching $X$ such that for all $n \in {\mathbb N}$,
$x_n=a$, $y_n=b$ where $a, b \in {\mathbb C}$ are fixed constants.
Let $d$ be the distance from $v$ to the edge $e_1\leftrightarrow
(X_1,Y_1)$. It follows that there is a node $v_n$ of ${\mathcal H}$
at distance $d$ from each of $e_n \leftrightarrow (X_n,Y_n)$. Since
$e_n \rightarrow X$, $v_n\rightarrow X$, and hence $X$ is the limit
point of distinct $\lambda_n \in \EE$. \qed

\begin{rmk}
We end this section by remarking that if $\rho \in {\mathcal X}^{\rm
disc}$ satisfies the conditions of Theorem \ref{thm:Efordiscrete},
and $\EE \neq \PL$, then the proof of Theorem \ref{thm:Efordiscrete}
implies that in fact, the stabilizer of $\rho$ in the mapping class
group $\Gamma$ of $T$ is relatively large, in particular, is not
finite or virtually cyclic. This answers a question posed by
Bowditch in \cite{bowditch1998plms}. Makoto Sakuma
\cite{sakuma_private} has also independently obtained examples of
characters whose stabilizer in $\Gamma$ is not finite or virtually
cyclic by considering the representations arising from the two
bridge knot complements; in these cases there are at least two
(hence infinitely many) $X \in \C$ for which $\phi(X)=0$.
\end{rmk}


\section{Proof of Theorem \ref{thm:EforImaginary}}\label{s:proofforimaginary} %
For the rest of this section, we fix $\rho \in {\mathcal
X}_{\kappa}^I$, where $\kappa < 2$, with corresponding Fricke trace
map $\phi$. It will be convenient and visually easier to regard
$\phi$ as a map from $\Omega(\Sigma)$ to $\mathbb C$. Recall from \S
\ref{ss:tricolor} that there is a partition of $\C$ and
$\Omega(\Sigma)$ into three equivalence classes $\R,\G$ and $\B$. We
use the letters $Z, Y$ and $X$  to denote the elements of $\R,\G$
and $\B$ respectively, and $W$ to denote a general element of $\C$.
Since $\rho$ is imaginary, $\phi$ takes purely imaginary values on
two of the equivalence classes and real values on the third; we may
assume that it takes real values on $Z \in \R$. For $X \in \B$, $Y
\in \G$ and $Z \in \R$, we use the convention $\phi(X)=ix$,
$\phi(Y)=iy$ and $\phi(Z)=z$, where $x,y,z \in {\mathbb R}$ for the
rest of this section (note that this is different from the
convention used earlier). Then if $(X,Y,Z) \in \GT$  and
$(X,Y;Z,Z'),(Y,Z;X,X'),(Z,X;Y,Y')\in \GQ$, the vertex relation
(\ref{eqn:vertexrelation}) and edge relation
(\ref{eqn:edgerelation}) can be rewritten as:
\begin{equation}\label{eqn:imaginaryvertex}
    -x^2-y^2+z^2+xyz-2=\kappa
\end{equation}
\begin{equation}\label{eqn:imaginaryedge}
    z+z'=-xy, \quad y+y'=xz, \quad x+x'=yz
\end{equation}

Note that  by our convention (\S \ref{sss:realimaginarycharacters})
$\rho \in {\mathcal X}_{\kappa}^I$ is not dihedral. Hence, by
Theorem \ref{thm:E=PL}, $\EE \neq \PL$, and in fact $\EE$ has empty
interior in $\PL$. We will be using the results of
Goldman-Stantchev, in particular, the classification result (Theorem
A of \cite{goldman-stantchev}) and the $\tau$-reduction algorithm
(\S 4 of \cite{goldman-stantchev}). The following strengthening of
Proposition \ref{prop:irrationalend} concerning irrational ends is
key.



\begin{figure}
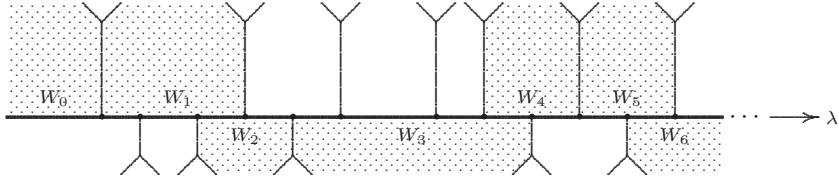

\begin{center}
\mbox{
\beginpicture
\setcoordinatesystem units <0.5in,0.5in> \setplotarea x from 1 to
10, y from 0 to 2.5 \plot 1 1 8.5 1 / \plot 1 0.99 8.5 0.99 /
\plot 1 1.01 8.5 1.01 / \put {\mbox{\Huge $\cdot$}} [cc]
<0mm,-0.1mm> at 2 1 \plot 2 1 2 2 / \plot 1.8 2.2 2 2 2.2 2.2 /
\put {\mbox{\Huge $\cdot$}} [cc] <0mm,-0.1mm> at 2.4 1 \plot 2.4 1
2.4 0.6 / \plot 2.2 0.4 2.4 0.6 2.6 0.4 / \put {\mbox{\Huge
$\cdot$}} [cc] <0mm,-0.1mm> at 3 1 \plot 3 1 3 0.6 / \plot 2.8 0.4
3 0.6 3.2 0.4 / \put {\mbox{\Huge $\cdot$}} [cc] <0mm,-0.1mm> at
3.5 1 \plot 3.5 1 3.5 2 / \plot 3.3 2.2 3.5 2 3.7 2.2 / \put
{\mbox{\Huge $\cdot$}} [cc] <0mm,-0.1mm> at 4 1 \plot 4 1 4 0.6 /
\plot 3.8 0.4 4 0.6 4.2 0.4 / \put {\mbox{\Huge $\cdot$}} [cc]
<0mm,-0.1mm> at 4.5 1 \plot 4.5 1 4.5 2 / \plot 4.3 2.2 4.5 2 4.7
2.2 / \put {\mbox{\Huge $\cdot$}} [cc] <0mm,-0.1mm> at 5.5 1 \plot
5.5 1 5.5 2 / \plot 5.3 2.2 5.5 2 5.7 2.2 / \put {\mbox{\Huge
$\cdot$}} [cc] <0mm,-0.1mm> at 6 1 \plot 6 1 6 2 / \plot 5.8 2.2 6
2 6.2 2.2 / \put {\mbox{\Huge $\cdot$}} [cc] <0mm,-0.1mm> at 6.5 1
\plot 6.5 1 6.5 0.6 / \plot 6.3 0.4 6.5 0.6 6.7 0.4 / \put
{\mbox{\Huge $\cdot$}} [cc] <0mm,-0.1mm> at 7 1 \plot 7 1 7 2 /
\plot 6.8 2.2 7 2 7.2 2.2 / \put {\mbox{\Huge $\cdot$}} [cc]
<0mm,-0.1mm> at 7.5 1 \plot 7.5 1 7.5 0.6 / \plot 7.3 0.4 7.5 0.6
7.7 0.4 / \put {\mbox{\Huge $\cdot$}} [cc] <0mm,-0.1mm> at 8 1
\plot 8 1 8 2 / \plot 7.8 2.2 8 2 8.2 2.2 / \put {\mbox{$\cdots$}}
[cc] <0mm,0mm> at 8.75 1 \arrow <6pt> [.16,.6] from 9 1 to 9.5 1
\setlinear \setshadegrid span <2pt> \vshade 1 1 2
<0pt,0pt,0pt,0pt> 3.5 1 2 / \setlinear \setshadegrid span <2pt>
\vshade 1 2 2.2 <0pt,0pt,0pt,0pt> 1.8 2 2.2 / \setlinear
\setshadegrid span <2pt> \vshade 1.8 2 2.2 <0pt,0pt,0pt,0pt> 2 2 2
/ \setlinear \setshadegrid span <2pt> \vshade 2 2 2
<0pt,0pt,0pt,0pt> 2.2 2 2.2 / \setlinear \setshadegrid span <2pt>
\vshade 3.3 2 2.2 <0pt,0pt,0pt,0pt> 3.5 2 2 / \setlinear
\setshadegrid span <2pt> \vshade 2.2 2 2.2 <0pt,0pt,0pt,0pt> 3.3 2
2.2 / \setlinear \setshadegrid span <2pt> \vshade 6 1 2
<0pt,0pt,0pt,0pt> 8 1 2 / \setlinear \setshadegrid span <2pt>
\hshade 2 7 8  <0pt,0pt,0pt,0pt> 2.2 7.2 7.8  / \setlinear
\setshadegrid span <2pt> \hshade 2 6 7  <0pt,0pt,0pt,0pt> 2.2 6.2
6.8  / \setlinear \setshadegrid span <2pt> \vshade 3 0.6 1
<0pt,0pt,0pt,0pt> 6.5 0.6 1  / \setlinear \setshadegrid span <2pt>
\hshade 0.4 3.2 3.8 <0pt,0pt,0pt,0pt> 0.6 3 4 / \setlinear
\setshadegrid span <2pt> \hshade 0.4 4.2 6.3 <0pt,0pt,0pt,0pt> 0.6
4 6.5 / \setlinear \setshadegrid span <2pt> \hshade 0.6 7.5 8.5
<0pt,0pt,0pt,0pt> 1 7.5 8.5 / \setlinear \setshadegrid span <2pt>
\hshade 0.4 7.7 8.5 <0pt,0pt,0pt,0pt> 0.6 7.5 8.5 /
\put {\mbox{\scriptsize $\lambda$}} [lc] <1mm,0mm> at 9.5 1 \put
{\mbox{\scriptsize $W_{0}$}} [cc] <0mm,0mm> at 1.5 1.2 \put
{\mbox{\scriptsize $W_{1}$}} [cc] <0mm,0mm> at 2.8 1.2 \put
{\mbox{\scriptsize $W_{2}$}} [cc] <0mm,0mm> at 3.5 0.8 \put
{\mbox{\scriptsize $W_{3}$}} [cc] <0mm,0mm> at 5.25 0.8 \put
{\mbox{\scriptsize $W_{4}$}} [cc] <0mm,0mm> at 6.5 1.2 \put
{\mbox{\scriptsize $W_{5}$}} [cc] <0mm,0mm> at 7.5 1.2 \put
{\mbox{\scriptsize $W_{6}$}} [cc] <0mm,0mm> at 8 0.8
\endpicture}
\end{center}

\caption{The sequence $\{W_{n}\}\longrightarrow \lambda$ along
$P$.}\label{fig:pathtoirrationalend}
\end{figure}


\begin{prop} \label{prop::irrationalendforimaginary}
Suppose $\rho \in {\mathcal X}_{\kappa}^I$, where $\kappa < 2$ and
$\lambda \in \EE$ with $\lambda \not\in \C$. Then there exists a
sequence of distinct elements $\{W_n\} \subset \C \cap \EE$ such
that $W_n \rightarrow \lambda$.
\end{prop}

\begin{pf}
By Proposition \ref{prop:irrationalend}, $\Omega(2) \neq \emptyset$.
Choose $W_0 \in \Omega(2)$ and let $P=\{e_n\}$ be the minimal path
in $\Sigma(\Omega)$ from $W_0$ limiting at $\lambda$. By Lemma
\ref{lem:quasiconvexity}, we can find a sequence $\{W_n\} \subset
\Omega(2)$  such that  for all $n \ge 1$, $W_n$ is adjacent to $P$,
$W_n,W_{n+1}$ are neighbors and  $W_n \rightarrow \lambda$ (see
Figure \ref{fig:pathtoirrationalend}). We consider two cases, first
when $W_n \in \R$ for infinitely many $n$, and
secondly, when $W_n \not\in \R$ for all but finitely many $n$.\\

{\it Case 1. $W_n \in \R$ for infinitely many $n$}. \, This gives a
subsequence $\{Z_n\} \subset \R$ approaching $\lambda$ such that
$|z_n| \le 2$. If there are infinitely many elements in this
sequence with $|z_n|<2$, these are in $\EE$ by Proposition
\ref{prop:rationalendinvariants}, so we are done. Hence we might as
well assume that $z_n=\pm 2$ for all $n$. Then there exists
infinitely many $k$ such that $\phi(W_k)=\pm 2$ in the sequence
$\{W_n\}$. Fix such a $k$ and consider the triple $W_{k-1}, W_k,
W_{k+1}$. Using the sign-change automorphism (\S
\ref{sss:signchange}) we might as well assume that $w_k=2$.
$W_{k-1}$ and $W_{k+1}$ are neighbors of $W_k$ (hence, not in $R$),
with $\phi(W_{k-1})=iw_{k-1}$, $\phi(W_{k+1})=iw_{k+1}$, where
$w_{k-1}, w_{k+1} \in [-2,2]$, by assumption. If $w_{k-1}=0$ (or
$w_{k+1}=0$), then $W_{k-1} \in \EE$ (resp. $W_{k+1} \in \EE$). If
both $w_{k-1}$ and $w_{k+1}$ have the same sign, then, since the
values of the neighbors of $W_k$ grow linearly (Lemma
\ref{lem:neighbors}(c)), the difference in the values of successive
neighbors of $W_k$ is $ci$, where $0<c<2$, and we can find two
successive neighbors, say $X$ and $Y$, of $W_k$ such that
$\phi(X)=ix$, $\phi(Y)=iy$ with $-2\le x<0<y\le 2$. If $w_{k-1}$ and
$w_{k+1}$ have different signs, then a simple argument using the
neighbors of $W_k$ in between $W_{k-1}$ and $W_{k+1}$ again yields
neighbors $X,Y$ of $W_k$ with the same property. If $x=-2$ and
$y=2$, we have $(w_k,ix,iy)=(2,-2i,2i)$, and $\rho$ can be easily
shown to satisfy the extended BQ-conditions (see
\cite{tan-wong-zhang2005gmm}), so that $\EE=\emptyset$,
contradicting the fact that $\rho$ has an irrational end. Hence, we
may assume that either $x>-2$ or $y<2$. Now if we consider the
generating quadruple $(X,Y;W_k,Z_k)$, we have $|z_k|=|-xy-2|<2$, so
we obtain $Z_k \in \R \cap \EE$. In all subcases,  we obtain a
subsequence $\{W_n\} \subset \C \cap \EE$ such that $W_n \rightarrow
\lambda$, which completes the proof in this case.

{\it Case 2. $W_n \not\in \R$ for all but finitely many $n$}. \,
Re-indexing if necessary, we might as well assume that $W_n \not\in
\R$ for all $n$, and  $\{W_n\}=X_1,Y_1,X_2,Y_2,\dots$, where $X_n
\in \B$, $Y_n \in \G$. The sequence $\{W_n\}$ must cross the path
$P$ infinitely often, since the end of $P$ is irrational, and $W_n
\not\in \R$ for all $n$. Hence, renaming $\G$ and $\B$ if necessary,
we obtain a subsequence of generating pairs $\{(X_n,Y_n)\}\subset
\GP$ and a nested sequence of closed intervals $[X_n,Y_n]$ such that
$$\lambda \in \cdots \subset [X_{n+1}, Y_{n+1}] \subset [X_n,Y_n]\cdots \subset [X_1,Y_1],$$
$X_n, Y_n \rightarrow \lambda$, and $|ix_n|, |iy_n| \le 2$ for all
$n \ge 1$. Passing again to a subsequence, we may assume that $x_n
\rightarrow a$, $y_n \rightarrow b$, where $a,b \in [-2,2]$. Let
$e_n \in E_r \subset E(\Sigma)$ be the edge of $P$ corresponding
to the generating pair $(X_n,Y_n)$ and $v_n \in V(\Sigma)$ be the
end of $e_n$ which is closer to $\lambda$. Note that the quantity
$\tau_n=-z_nz_n'$ associated to the edge $e_n \leftrightarrow
(X_n,Y_n;Z_n,Z_n') \in \GQ$ approaches a fixed constant $c$ which
depends only on $a$ and $b$. We now consider the case where both
$a,b \neq 0$ and the case where one of $a,b=0$ separately, and
perform the $\tau$-reduction algorithm at each $v_n$.

\emph{Subcase} (i): $a, b \neq 0$. \,Then, starting at $v_n$, the
$\tau$-reduction algorithm of Goldman-Stantchev
\cite{goldman-stantchev}, see also the Appendix, produces a sequence
in $\C$ terminating after a finite number of steps at some $Z_n \in
\EE$. For $n$ sufficiently large, since $a,b \neq 0$, the first
(two) steps of the algorithm reduces $\tau$ by at least some
positive constant depending only on $a$ and $b$, hence if $n$ is
sufficiently large, the algorithm starting at $v_n$ does not cross
the edges $e_{n-1}\leftrightarrow (X_{n-1}, Y_{n-1})$ and
$e_{n+1}\leftrightarrow (X_{n+1}, Y_{n+1})$. Hence, $Z_n \in
[X_{n-1}, Y_{n-1}]\setminus [X_{n+1}, Y_{n+1}]$. Again passing to a
subsequence, we obtain a sequence of distinct $Z_n \in \EE$
approaching $\lambda$.

\emph{Subcase} (ii): One of $a,b =0$, say $a=0$. \,Then for any
$\epsilon >0$, there exists $N$ such that $|ix_n| <\epsilon$ for all
$n>N$. Again, the $\tau$-reduction algorithm starting at $v_n$
terminates, after a finite number of steps, at some $Z_n \in \C$
such that $Z_n \in \EE$. Furthermore, if $|ix_n|$ is sufficiently
small, by \cite{goldman-stantchev} (see Appendix), $Z_n$ is in fact
a neighbor of $X_n$. Again it follows that for $n$ sufficiently
large, $Z_n\in [X_{n-1}, Y_{n-1}]\setminus [X_{n+1}, Y_{n+1}]$, and
the conclusion follows as in Subcase (i).
%
%
%
%
%
\end{pf}

\begin{figure}
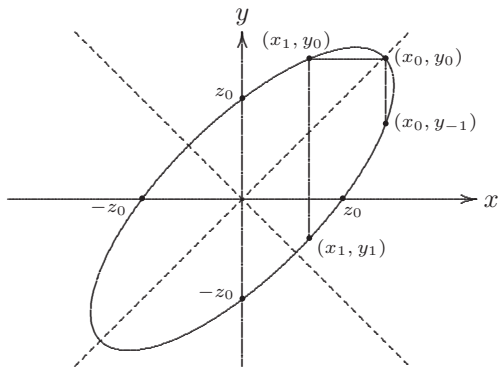


\begin{center}
\mbox{
\beginpicture
\setcoordinatesystem units <0.35in,0.35in> \setplotarea x from -4
to 4, y from -3.5 to 3.5 \arrow <6pt> [.16,.6] from 0 -2.5 to 0
2.5 \put {$y$} [cb] <0mm,1mm> at 0 2.5 \arrow <6pt> [.16,.6] from
-3.5 0 to 3.5 0 \put {$x$} [lt] <1mm,0.5mm> at 3.5 0
\startrotation by .7071067810 .7071067810 about  0 0
\ellipticalarc axes ratio 1:.3779644731 360 degrees from 3 0
center at 0 0 \stoprotation \setdashes<2pt> \plot -2.5 -2.5 2.5
2.5 / \plot -2.5 2.5 2.5 -2.5 / \setsolid \put {\mbox{\Huge
$\cdot$}} [cc] <0mm,-0.1mm> at 1 -.5962912020 \put {\mbox{\Huge
$\cdot$}} [cc] <0mm,-0.1mm> at 1 2.096291202 \put {\mbox{\Huge
$\cdot$}} [cc] <0mm,-0.1mm> at 2.144436802 2.096291202 \put
{\mbox{\Huge $\cdot$}} [cc] <0mm,-0.1mm> at 2.144436802
1.120364001 \plot 1 -.5962912020 1 2.096291202 2.144436802
2.096291202 2.144436802 1.120364001 / \put {\mbox{\Huge $\cdot$}}
[cc] <0mm,-0.1mm> at 0 1.5 \put {\mbox{\scriptsize $z_0$}} [rb]
<-1mm,0mm> at 0 1.5 \put {\mbox{\Huge $\cdot$}} [cc] <0mm,-0.1mm>
at 0 -1.5 \put {\mbox{\scriptsize $-z_0$}} [rb] <-1mm,0mm> at 0
-1.5 \put {\mbox{\Huge $\cdot$}} [cc] <0mm,-0.1mm> at 1.5 0 \put
{\mbox{\scriptsize $z_0$}} [lt] <0mm,-1mm> at 1.5 0 \put
{\mbox{\Huge $\cdot$}} [cc] <0mm,-0.1mm> at -1.5 0 \put
{\mbox{\scriptsize $-z_0$}} [rt] <-2mm,-0.2mm> at -1.5 0 \put
{\mbox{\scriptsize $(x_1,y_1)$}} [lt] <1mm,0mm> at 1 -.5962912020
\put {\mbox{\scriptsize $(x_1,y_0)$}} [cb] <-1.7mm,1.2mm> at 1
2.096291202 \put {\mbox{\scriptsize $(x_0,y_0)$}} [lc] <1mm,0.2mm>
at 2.144436802 2.096291202 \put {\mbox{\scriptsize
$(x_0,y_{-1})$}} [lc] <1mm,0.2mm> at 2.144436802 1.120364001

\endpicture}

\end{center}
\caption{The ellipse $x^2-z_0xy+y^2=z_0^2$ ($z_0=1.5$) and
successive values of $x_n, y_n$.}\label{fig:ellipse}
\end{figure}

{\sl Proof of Theorem \ref{thm:EforImaginary}.} \,The proof for
Theorem \ref{thm:EforImaginary} now follows easily, along similar
lines to the proof of Theorem \ref{thm:EforReal}. If $|\EE \cap \C|
>1$, it follows that $|\EE \cap \C|=\infty$ by Lemma
\ref{lem:neighbors}(b), and furthermore, each $W \in \EE \cap \C$
is an accumulation  point of $\EE$. If $\EE$ has an irrational end
$\lambda$, it follows from Proposition
\ref{prop::irrationalendforimaginary} that $|\EE \cap \C|=\infty$
and $\lambda$ is an accumulation point of $\EE$. Hence, $\EE$ is
either empty, has one element $W \in \C$ or is a Cantor set. The
fact that $\EE=\emptyset$ occurs only for $\kappa \le -14$ and not
for $-14<\kappa <2$ follows from \cite{goldman-stantchev}
(actually, the case $\kappa=-14$ was not covered, but as we saw
earlier, the example $\rho \in \X_{-14}$ where $\iota (\rho)=(2,
2i,-2i)$ has $\EE=\emptyset$). It is fairly easy to construct, for
all $\kappa <2$, examples where $\EE =\{X\}$: we use $\rho$ with
$\iota(\rho)=(0,yi,z)$ where $y,z$ are chosen such that $y,z>2$
and $z^2-y^2=\kappa+2$. In this case $\C(2)=\EE=\{X\}$. Similarly,
it is easy to construct examples where $|\EE|>1$, and hence is a
Cantor set: we just need to make sure that there are at least two
element in $\C \cap \EE$. We leave this to the reader. Finally, to
prove part (i) of the Theorem, we first show that for $\rho \in
\X_{-2}$, if $\EE=\{W\}$, then $W \not\in \R$. Suppose not, then
we have some $(X_0,Y_0,Z_0) \in \GT$ such that $z_0 \in (-2,2)$.
If $z_0=0$, then $ix_0=iy_0=0$, a contradiction since $\EE$ has
only one element. Hence, using the sign change automorphism if
necessary, we may assume that $z_0 \in (0,2)$. Use $X_n,Y_n$ to
denote the successive neighbors of $Z_0$, with values $ix_n$,
$iy_n$ respectively. Then $x_n,y_n$ lie on the ellipse
$$x^2-z_0xy+y^2=z_0^2,$$ with major axis $y=x$ and minor axis
$y=-x$, and intercepts $(\pm z_0,0)$ and $(0, \pm z_0$). The values
$x_n,y_n$ are obtained by starting at the point $(x_0,y_0)$ on the
ellipse and taking the coordinates of the  successive intersections
of the ellipse with the up/down and left/right path; see Figure
\ref{fig:ellipse} (compare with \cite{goldman-stantchev}). There is
at least one intercept which lies in either the second or the fourth
quadrant, that is, there exist some successive neighbors $X_n,Y_n$
(or $Y_n,X_{n+1}$) of $Z_0$ such that $x_n, y_n$ have opposite signs
and $|x_n|,|y_n|\le z_0<2$. However, in this case, for the
generating quadruple $(X_n,Y_n; Z_0, Z_n) \in \GQ$,
$|z_n|=|x_ny_n-z_0|<2$, contradicting the fact that $\EE$ has only
one element. Hence we conclude that the single element $W$ of $\EE$
is not in $\R$; so it must be in $\G$ or $\B$, which implies $w=0$.
The conclusion then follows easily from the vertex relation
(\ref{eqn:vertexrelation}) and the connectedness of $\C(2)$. \qed



%
%
%
%


\vskip 15pt
\section{Appendix: The $\tau$-reduction algorithm}\label{s:tracereduction}%
\vskip 5pt

We give a brief description of the $\tau$-reduction algorithm given
by Goldman and Stantchev in \cite{goldman-stantchev}. Bowditch's
algorithm in \cite{bowditch1998plms} for type-preserving imaginary
Markoff maps ($\rho \in {\mathcal X}_{-2}^I$) is essentially a
special case of this, although it was couched in a different setting
using the Jorgenson parameters, $a=x/yz, b=y/xz$ and $c=z/xy$
instead of $x,y,z$. We will describe the results and algorithm but
refer the reader to \cite{goldman-stantchev} for detailed proofs.

Fix a $\rho \in {\mathcal X}_{\kappa}^I$, where $\kappa <2$, with
corresponding $\phi$. We adopt the notation in \S
\ref{s:proofforimaginary}; in particular, we use $X$, $Y$ and $Z$
for the elements of $\B, \G$ and $\R$ respectively (thought of as
subsets of $\Omega(\Sigma)$), and write $z=\phi(Z)$, $ix=\phi(X)$,
$iy=\phi(Y)$ where $x,y,z \in {\mathbb R}$. If $e \in E_r \subset
E(\Sigma)$ corresponds to the generating quadruple $(X,Y;Z,Z')$,
define $\tau(e)=-zz'$. Similarly, for $v \in V(\Sigma)$, define
$\tau(v)=\tau(e)$ where $e$ is the unique edge in $E_r$ with one
endpoint at $v$. Note that $\tau$ is invariant under the sign-change
automorphism.

The $\tau$-reduction algorithm has as a starting point a vertex
$v_0\in V(\Sigma)$, and produces a (unique) finite sequence of
adjacent vertices $v_0, v_1, \cdots, v_N$ such that either
\begin{enumerate}
\item[(a)] $\tau(v_N)$ in minimal among all vertices $v \in
V(\Sigma)$, $v_N$ is an attractor for the flow $f_{\rho}$, $\phi$
satisfies the extended BQ-conditions and $\EE=\emptyset$; or

\item[(b)] one of the three regions adjacent to $v_N$ is in $\EE$.

\end{enumerate}




\begin{figure}
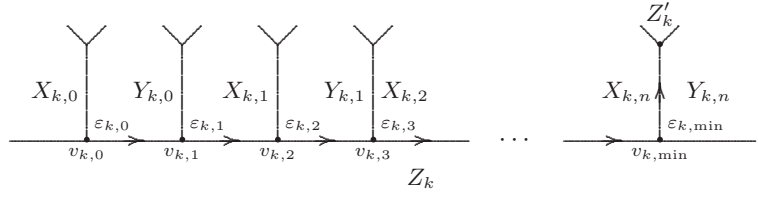

\begin{center}
\mbox{
\beginpicture
\setcoordinatesystem units <0.5in,0.5in> \setplotarea x from -1 to
7, y from -0.5 to 1.5 \plot -0.8 0 4 0 / \plot 0 0 0 1 / \plot -0.2
1.2 0 1 0.2 1.2 / \put {\mbox{\Huge $\cdot$}} [cc] <0mm,-0.1mm> at 0
0 \put {\mbox{\scriptsize $v_{k,0}$}} [ct] <0mm,-1mm> at 0 0 \put
{\mbox{\small $X_{k,0}$}} [rb] <-1mm,1mm> at 0 0.3 \put
{\mbox{\scriptsize $\varepsilon_{k,0}$}} [lb] <1mm,1mm> at 0 0
\arrow <6pt> [.16,.6] from 0.4 0 to 0.6 0
\plot 1 0 1 1 / \plot 0.8 1.2 1 1 1.2 1.2 / \put {\mbox{\Huge
$\cdot$}} [cc] <0mm,-0.1mm> at 1 0 \put {\mbox{\scriptsize
$v_{k,1}$}} [ct] <0mm,-1mm> at 1 0 \put {\mbox{\small $Y_{k,0}$}}
[rb] <-1mm,1mm> at 1 0.3 \put {\mbox{\scriptsize
$\varepsilon_{k,1}$}} [lb] <1mm,1mm> at 1 0 \arrow <6pt> [.16,.6]
from 1.4 0 to 1.6 0
\plot 2 0 2 1 / \plot 1.8 1.2 2 1 2.2 1.2 / \put {\mbox{\Huge
$\cdot$}} [cc] <0mm,-0.1mm> at 2 0 \put {\mbox{\scriptsize
$v_{k,2}$}} [ct] <0mm,-1mm> at 2 0 \put {\mbox{\small $X_{k,1}$}}
[rb] <-1mm,1mm> at 2 0.3 \put {\mbox{\scriptsize
$\varepsilon_{k,2}$}} [lb] <1mm,1mm> at 2 0 \arrow <6pt> [.16,.6]
from 2.4 0 to 2.6 0
\plot 3 0 3 1 / \plot 2.8 1.2 3 1 3.2 1.2 / \put {\mbox{\Huge
$\cdot$}} [cc] <0mm,-0.1mm> at 3 0 \put {\mbox{\scriptsize
$v_{k,3}$}} [ct] <0mm,-1mm> at 3 0 \put {\mbox{\small $Y_{k,1}$}}
[rb] <-1mm,1mm> at 3 0.3 \put {\mbox{\scriptsize
$\varepsilon_{k,3}$}} [lb] <1mm,1mm> at 3 0 \arrow <6pt> [.16,.6]
from 3.4 0 to 3.6 0 \put {\mbox{\small $X_{k,2}$}} [lb] <1mm,1mm>
at 3 0.3
\put {\mbox{$\cdots$}} [cc] <0mm,0mm> at 4.5 0 \plot 5 0 7 0 /
\arrow <6pt> [.16,.6] from 5.4 0 to 5.6 0 \plot 6 0 6 1 / \plot
5.8 1.2 6 1 6.2 1.2 / \put {\mbox{\Huge $\cdot$}} [cc]
<0mm,-0.1mm> at 6 0 \put {\mbox{\scriptsize $v_{k,{\rm min}}$}}
[ct] <0mm,-1mm> at 6 0 \put {\mbox{\small $X_{k,n}$}} [rb]
<-1mm,1mm> at 6 0.3 \put {\mbox{\scriptsize $\varepsilon_{k,{\rm
min}}$}} [lb] <1mm,1mm> at 6 0 \put {\mbox{\small $Y_{k,n}$}} [lb]
<1mm,1mm> at 6.2 0.3 \arrow <6pt> [.16,.6] from 6 0.4 to 6 0.6
\put {\mbox{\Huge $\cdot$}} [cc] <0mm,-0.1mm> at 6 1 \put
{\mbox{\small $Z'_{k}$}} [cb] <0mm,1mm> at 6 1.12 \put
{\mbox{\small $Z_{k}$}} [cb] <0mm,0mm> at 3.5 -0.5
\endpicture}\end{center}

\caption{$\tau$-algorithm moves along $\partial
Z_{k}$}\label{fig:neigborsofZ}
\end{figure}

This also produces a sequence of elements $W_1,W_2, \cdots,W_N \in
\Omega(\Sigma)$, where $W_k$ is the (unique) region adjacent to
$v_k$ which is not adjacent to $v_{k-1}$. The algorithm consists of
two general types of moves, the first type moves the vertices around
the boundary of $Z \in \R$ until we reach a point where $\tau$ is
minimum among all the vertices lying on the boundary of $Z$, then
the second type ``flips'' the vertex across the edge $e \in E_r$ to
a vertex which is now adjacent to a different $Z' \in \R$. The main
point is that the algorithm is uniquely determined, does not
backtrack, and terminates after a finite number of steps in one of
the two possibilities (a) or (b) listed above.

\vskip 3pt

\noindent\emph{The algorithm:} Let $v_0 \leftrightarrow
(X_0,Y_0,Z_0) \in \GT$, and $e\leftrightarrow (X_0,Y_0;Z_0,Z_0') \in
\GQ$, where $e \in E_r$ is adjacent to $v_0$. If $ix_0, iy_0$ or
$z_0 \in (-2,2)$, then $N=0$ and we are done (case (b) above). If
$|z|, |z'| \ge 2$ and both $z_0$, $z_0'$ have the same signs then
$N=0$ if $|z_0|\le|z_0'|$, and $N=1$ if $|z_0|>|z_0'|$, where $v_1
\leftrightarrow (X_0,Y_0,Z_0')$. By the results of
\cite{goldman-stantchev} and \cite{tan-wong-zhang2005gmm}, $v_N$ is
a sink which is an attractor for the flow $f_{\rho}$ (see also the
properties listed below), and we are done (case (a) above). If none
of the above is satisfied, we perform the inductive step below.

\vskip 3pt

\begin{figure}
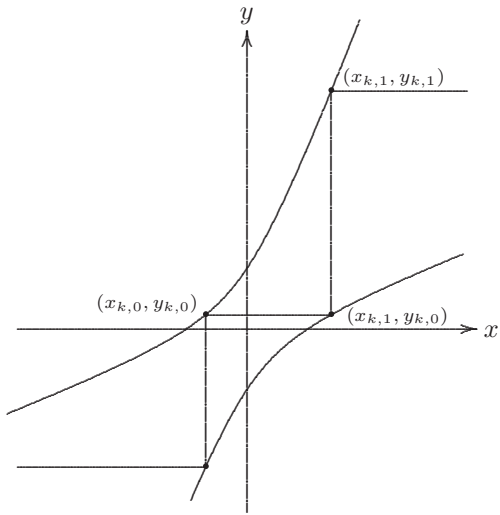


\begin{center}
\mbox{
\beginpicture
\setcoordinatesystem units <0.12in,0.12in> \setplotarea x from -10
to 10, y from -8 to 14 \arrow <6pt> [.16,.6] from 0 -8 to 0 13
\put {$y$} [cb] <0mm,1mm> at 0 13 \arrow <6pt> [.16,.6] from -10 0
to 10 0 \put {$x$} [lt] <1mm,0.5mm> at 10 0 \startrotation by
.7071067810 .7071067810 about  0 0 \setquadratic \plot -10.
4.774934554 -8.850000000 4.297033861 -7.700000000 3.828576760
-6.550000000 3.373499667 -5.400000000 2.938026548 -4.250000000
2.532291453 -3.100000000 2.173016336 -1.950000000 1.886928721
-.8000000000 1.711139971 .3500000000 1.680624884 1.500000000
1.802775638 2.650000000 2.050487747 3.800000000 2.384952829
4.950000000 2.774977477 6.100000000 3.200312485 7.250000000
3.648629880 8.400000000 4.112420212 9.550000000 4.586992479
10.70000000 5.069319480 11.85000000 5.557382477 13. 6.049793384 /
\plot -7. -3.549647870 -6.200000000 -3.238518180 -5.400000000
-2.938026548 -4.600000000 -2.651791847 -3.800000000 -2.384952829
-3. -2.144761058 -2.200000000 -1.941133689 -1.400000000
-1.786616915 -.6000000000 -1.694697613 .2000000000 -1.675708805 1.
-1.732050808 1.800000000 -1.856879102 2.600000000 -2.037645700
3.400000000 -2.260973242 4.200000000 -2.515551629 5. -2.792848008
5.800000000 -3.086745859 6.600000000 -3.392933834 7.400000000
-3.708368914 8.200000000 -4.030880796 9. -4.358898944 /
\stoprotation \setlinear \put {\mbox{\Huge $\cdot$}} [cc]
<0mm,-0.1mm> at -1.8 -6.024154028 \put {\mbox{\Huge $\cdot$}} [cc]
<0mm,-0.1mm> at -1.8 0.6241540277 \put {\mbox{\Huge $\cdot$}} [cc]
<0mm,-0.1mm> at 3.672462083 0.6241540277 \put {\mbox{\Huge
$\cdot$}} [cc] <0mm,-0.1mm> at 3.672462083 10.39323222 \plot -10
-6.024154028 -1.8 -6.024154028 -1.8 0.6241540277 3.672462083
0.6241540277 3.672462083 10.39323222 9.6 10.39323222 / \put
{\mbox{\scriptsize $(x_{k,0},y_{k,0})$}} [rb] <-1mm,0.2mm> at -1.8
0.6241540277 \put {\mbox{\scriptsize $(x_{k,1},y_{k,0})$}} [lc]
<2mm,-0.2mm> at 3.672462083 0.6241540277 \put {\mbox{\scriptsize
$(x_{k,1},y_{k,1})$}} [lb] <1.5mm,0.5mm> at 3.672462083
10.39323222
\endpicture}
\end{center}

\caption{The hyperbola $-x^2-y^2+z_k(xy)+z_k^2=\kappa+2$
($z_k=3,\kappa=0$) and successive values of $x_{k,n},
y_{k,n}$.}\label{fig:zigzag}
\end{figure}

\noindent \emph{Inductive step:} Fix $k $ and suppose $v_k
\leftrightarrow (X_k,Y_k,Z_k)$ and $|z_k| \ge 2$. Let
$X_{k,n},Y_{k,n}$, $n \in {\mathbb Z}$ be the successive neighbors
of $Z_k$, where $X_k=X_{k,0}$ and $Y_k=Y_{k,0}$. Let $v_{k, 2n}
\leftrightarrow (Z_k, X_{k,n},Y_{k,n}) \in \GT$ and $v_{k,2n+1} \in
V(\Sigma) \leftrightarrow (Z_k, Y_{k,n},X_{k,n+1})$ so that
$\{v_{k,n}\}_{n \in {\mathbb Z}}$ are the successive vertices along
$\partial Z_k$. Let $\varepsilon_{k,n}\in E_r $ be the edge with
endpoint at $v_{k,n}$ which is not adjacent to $Z_k$ (see Figure
\ref{fig:neigborsofZ}). Using the sign change automorphism, we may
assume that $z_k \ge 2$. Then $(x_{k,n},y_{k,n})$ and
$(x_{k,n},y_{k,n-1})$ are coordinates of points on the hyperbola
(two parallel lines if $z_k=2$)
\begin{equation}\label{eqn:hyperbola}
-x^2-y^2+z_k(xy)+z_k^2=\kappa+2.
\end{equation}
The successive values of $x_{k,n},y_{k,n}$ in either direction can
be obtained by looking at the intercepts of the up/down left/right
zigzag path with the hyperbola (\ref{eqn:hyperbola}), starting from
the point $(x_0,y_0)=(x_{k,0}, y_{k,0})$ (see Figure
\ref{fig:zigzag}). We have the following facts:
\begin{enumerate}

\item [(i)] $x_{k,n}, y_{k,n} \rightarrow \pm \infty$ as $n
\rightarrow \pm \infty$.

\item [(ii)] Either exactly  one of $x_{k,n}$, $y_{k,n}=0$, or
there is a unique vertex $v_{k, min}$ adjacent to $Z_k$ at which
there is a change of signs in successive values of $x_{k,n},
y_{k,n}$, that is, the other two regions adjacent to $v_{k, min}$
have purely imaginary values with different signs.

\item [(iii)] If $x_{k,n}, y_{k,n} \neq 0$ for all $n$, then
$\tau(v_{k,n})>\tau(v_{k,min})$ with strict inequality, for all
$v_{k,n} \neq v_{k,min}$. The same holds if exactly  one of
$x_{k,n}$, $y_{k,n}=0$, except now $\tau$ achieves its minimum on
the two vertices on $\partial Z_k$ which are adjacent to the
region with value 0.

\item [(iv)] If $\{e_{k,n}\}_{n \in {\mathbb Z}}$ are the
successive edges of $\partial Z_k$, then $f_{\rho}(e_{k,n})$ is
directed towards $v_{k, min}$ for all $n \in {\mathbb Z}$, or the
region $X_{k,l}$ (or $Y_{k,l}$) neighboring $Z_k$ with value 0 (in
this case there is ambiguity in the definition of $f_{\rho}$ for
the edge between $Z_k$ and the region $X_{k,l}$ (or $Y_{k,l}$)
with value 0). This also follows from Proposition
\ref{prop:twoarrowspointingout}.
\end{enumerate}



Now starting at $v_k$, define a sequence of successive vertices
along $\partial Z_k$ until by (ii) above, we either reach a vertex
which has an adjacent region with value zero, or we reach $v_{k,
min}$. In the first case, we end the algorithm at that point as we
are in case (b). In the second case, if $(X_{k,n}, Y_{k,n}, Z_k,
Z_k') \leftrightarrow \varepsilon_{k,min}$, we check the value of
$z_k'=-x_{k,n}y_{k,n}-z_k$. Note that $x_{k,n}$ and $y_{k,n}$ have
different signs, so $-x_{k,n}y_{k,n}>0$, and that $z_k>2$ by
assumption. If $z_k'\ge z_k \ge 2$, we stop, as we are then in case
(a); otherwise, we perform the second type of move and the next
vertex is $v_{k,min}'$ which is adjacent to $v_{k,min}$ along
$\varepsilon_{k,min}$. If $z_k'\ge 2$, we stop as again, we are in
case (a). Similarly, if $z_k' \in(-2,2)$, we stop, since we are in
case (b). Finally, if $z_k' \le -2$, then $|z_k'|<|z_k|$. We now
proceed inductively with the new vertex $v_{k,min}'$ taking the
place of $v_k$, where we first perform  a sign-change automorphism
so we are back to the situation of $v_k$. It is clear that the
algorithm does not backtrack. To show that it terminates after a
number of steps, we need to consider two cases. If $|X|, |Y|$ is
bounded away from 0 for all $X \in \B$ and $Y \in \G$, then for each
step of the algorithm which is not a move along an edge of $E_r$,
$\tau$ reduces by at least a fixed amount $K>0$ so the algorithm
must terminate in a finite number of steps. Otherwise, there exists
$\epsilon >0$ (depending only on $\kappa$) such that if
$v_k\leftrightarrow (X_k,Y_k,Z_k)$ and $0<|ix_k|<\epsilon$, then
$v_{k+1}, \cdots, v_N$ all lie on $\partial X_k$, and the algorithm
terminates after a finite number of steps to case (b).

\vskip 30pt

{}

\end{document}